\documentclass[a4paper,10pt]{amsart}

\newcommand{\fx}{\mbox{\boldmath  $x$}}
\newcommand{\fw}{\mbox{\boldmath  $w$}}
\newcommand{\fv}{\mbox{\boldmath  $v$}}
\newcommand{\fs}{\mbox{\boldmath  $s$}}
\newcommand{\fe}{\mbox{\boldmath  $e$}}
\newcommand{\fu}{\mbox{\boldmath  $u$}}
\newcommand{\fU}{\mbox{\boldmath  $U$}}
\newcommand{\fb}{\mbox{\boldmath  $b$}}
\newcommand{\fzero}{\mbox{\boldmath  $0$}}

\usepackage{tikz}
\usepackage{cancel}
\usepackage{pgfplots}
\usepackage{cleveref}
\pgfplotsset{compat=newest}
\usetikzlibrary{plotmarks}

\usepackage{graphicx}
\usepackage{color}
\usepackage{xcolor}
\usepackage{amsmath,mathtools}
\usepackage{amsfonts}
\usepackage{mathrsfs}
\usepackage{color}
\usepackage{enumerate}
\usepackage{algorithm}
\usepackage{mdframed}
\usepackage{algorithmic}
\usepackage{siunitx}
\usepackage{multirow}
\usepackage{stackrel}
\usepackage{booktabs}
\usepackage{tabularx}
\usepackage{arydshln}
\usepackage{cite}
\usepackage{url}

\definecolor{darkgreen}{rgb}{.39 0.84 0.10}

\newcommand{\bgamma}{\boldsymbol{\gamma}}

\usepackage{changes}
\definechangesauthor[color=orange]{AP}

\newtheorem{Theorem}{Theorem}

\newtheorem{Corollary}{Corollary}
\newtheorem{Definition}{Definition}

\newtheorem{Lemma}{Lemma}

\newtheorem{Remark}{Remark}

\definecolor{g}{RGB}{1, 100, 42}
\newcommand{\R}{\mathbb{R}}

\usepackage[top=3.0cm,bottom=3.0cm,left=2.2cm,right=2.2cm]{geometry}
\usepackage{lineno}

\begin{document}
\title[Eigenvalues of multiple saddle-point matrices with block-triangular preconditioners]{Eigenvalue bounds for preconditioned symmetric multiple saddle-point matrices with block-triangular preconditioners}
\author{Luca Bergamaschi \and Michele Bergamaschi \and John W. Pearson
}
\footnotetext[1]{Department of Civil Environmental and Architectural Engineering, University of Padua, Via Marzolo, 9, 35100 Padua, Italy,
E-mail: \texttt{luca.bergamaschi@unipd.it}}
\footnotetext[2]{E-mail: \texttt{michele15456@gmail.com}}
\footnotetext[3]{School of Mathematics, The University of Edinburgh, James Clerk Maxwell Building, The King's Buildings, Peter Guthrie Tait Road, Edinburgh, EH9 3FD, United Kingdom,
E-mail: \texttt{j.pearson@ed.ac.uk}}


\begin{abstract}
We develop eigenvalue bounds for symmetric, block-tridiagonal multiple saddle-point linear systems, preconditioned with block-triangular matrices,
	based on approximate Schur complements. Irrespective on the number of blocks, we prove that all complex eigenvalues, with nontrivial imaginary part, 
    are strictly contained in a circle within the complex plane, with center $1$. 
    The real and positive eigenvalues are bounded in terms of the extremal roots of a sequences of parametric polynomials.
    Numerical results reveal that the bounds describe very well the eigenvalue  distribution of the preconditioned matrix.
\end{abstract}
\maketitle

\bigskip


\section{Introduction}
\label{sec:introduction}

We consider the iterative solution of a block-tridiagonal multiple saddle-point linear system
\begin{equation}\label{Eq1FULL}
\mathcal {A} \fx = \fb, \qquad \text{where} \quad
        \mathcal{A} = \begin{bmatrix}
        A_0 & B_1^\top &   &        &    \\
        B_1 & -A_1 & B_2^\top &        &        \\
          &  B_2 &  A_2 &          \ddots &\\
                &      & \ddots          &  \ddots       & B_N^\top \\
                &        &   & B_N &   (-1)^N A_N
\end{bmatrix}.
\end{equation}
We assume that $A_0\in\mathbb{R}^{n_0\times n_0}$ is symmetric positive definite, all other square block matrices $A_k\in\mathbb{R}^{n_k\times n_k}$ are symmetric positive semi-definite, and $B_k\in\mathbb{R}^{n_k\times n_{k-1}}$ have full rank ({with $n_k \le n_{k-1}$,} for $k=1,\dots, N$). 
{Denoting with  
$S_1 = A_1 + B_1 A_0^{-1} B_1^\top$ and $S_k = A_k + B_k{S}^{-1}_{k-1}B_k^\top$ (for $k = 2, \ldots, N$) the Schur complements, 
these conditions ensure that $S_k$
(for $k = 1, \ldots, N$) are invertible,
which also guarantees invertibility of $\mathcal{A}$ \cite{Bradley,BEIK2024403,pearson2023symmetric}. }

Linear systems involving the matrix $\mathcal{A}$ arise with $N =2$ (double saddle-point systems) in many scientific applications including magma--mantle dynamics \cite{Rhebergen},
liquid crystal director modeling \cite{RamGar2023}, or in the coupled Stokes--Darcy problem \cite{greifhe2023, Szyld, BeikBenzi2022, greif2026}, and the preconditioning of such linear systems has been considered in \cite{Balani-et-al-2023a, Balani-et-al-2023b, Benzi2018}, for instance.
In particular, block-diagonal preconditioners for the matrix {$\mathcal{A}$} have been thoroughly examined in
\cite{FRIGO202136,SZ,Bradley,PPNLAA24,bergamaschi2026eigenvalue,Mardal2026-MMMAS}. 
Augmented Lagrangian preconditioners have also been investigated in \cite{Benzi2026} to solve linear systems arising from
finite element discretizations of the fictitious domain method. 

Multiple saddle-point linear systems with $N > 2$ have recently attracted the attention of a number
of researchers. Such systems often stem from modeling multiphysics processes, i.e., the simultaneous simulation of different aspects of
physical systems and the interactions among them.
{
A recent example of this is the so-called four-field Thermo-Poroelasticity model which gives rise to a $4 \times 4$ block linear system.
Discretization and preconditioning issues are considered, e.g., in \cite{Mardal2026}.}
Other examples of multiple saddle-point linear systems can be found in \cite{BSZ2020,FerFraJanCasTch19,BerFerMar26}.
%

The accurate and efficient solution of linear systems in which $\mathcal A$ takes the form \eqref{Eq1FULL} is required to
obtain a reliable numerical simulation of these applications. Due to its size and sparsity, direct (factorization) methods
such as the $LDL^\top$ factorization of $\mathcal A$ (see, e.g., \cite{Greif2017} for more details)
are not recommended, since they can 
produce an excessively dense  triangular factor which cannot easily be stored when the size of the problems is very large. 

A common approach to overcome this issue is to consider iterative methods for the solution of $\mathcal A \fx = \fb$.
Among this class of methods, the most popular fall within the class of Krylov subspace solvers \cite{saadbook}, which exploit sparsity.
The matrix $\mathcal A$ in \eqref{Eq1FULL} is indefinite, meaning it has both negative and positive eigenvalues, posing challenges for the convergence of such iterative methods; see,
e.g., \cite{BenziGolubLiesen2005} and the references therein. If the intervals of negative and positive eigenvalues have large length, and/or the eigenvalues of smallest magnitude are close to zero, these features are expected to be highly detrimental to the numerical convergence. 
For this reason, preconditioning of block systems is invaluable, through devising a suitable preconditioner, $\mathcal P$ and transforming the original system into an equivalent system $\mathcal P^{-1} \mathcal A \fx = \mathcal P^{-1} \fb$ or $\mathcal A \mathcal P^{-1} (\mathcal P \fx) = \fb$, ideally with more favourable spectral properties.


To efficiently precondition $\mathcal A$, a number of different strategies can be followed, which take into account the block structure of the matrix. One class of approaches involves finding a symmetric positive definite preconditioner $\mathcal P$ of $\mathcal A$, with the aim of using short-term recurrence
iterative solvers, such as the MINRES \cite{minres} method, the prototypical Krylov subspace methods for such settings. 
In this case, the eigenvalues of the preconditioned system again lie in two intervals, one negative and one positive, with the action of an effective preconditioner improving the distribution of these intervals. Another strategy seeks a {symmetric indefinite or non-symmetric} preconditioner to be employed within, e.g., the GMRES method \cite{saad1986gmres}, with a block structure and spectral properties that mimic those of $A$. In this case, the eigenvalues of the preconditioned system can be clustered around $1$, within a subset of the complex right half-plane. 

Examples of the first class include block-diagonal preconditioners \cite{BenziGolubLiesen2005,desturler,bergamaschi2026eigenvalue,Pilotto-et-al2026} and a recent approach \cite{pearson2023symmetric,BMPP_COAP25,BB2026}
which, based on a decomposition $\mathcal A = \mathcal L\mathcal D \mathcal L^\top$ of multiple saddle-point systems, sets the preconditioner as (an approximation of)
$\mathcal  P = \mathcal L |\mathcal D| \mathcal L^\top$.
The block-triangular preconditioner that we analyze in this work for multiple saddle-point linear systems, belongs to the second class.  This preconditioner has been studied in \cite{Sim04}
for saddle-point linear systems, for which the eigenvalue distribution of the preconditioned matrix is described. Block-triangular preconditioners are also commonly used
to accelerate iterative solvers for double saddle-point linear systems; see, e.g., \cite{Benzi2018}.


The novelty of this work consists of extending the analysis of the block-triangular preconditioning approach for double saddle-point systems in \cite{Balani-et-al-2023b,BerFerMar26} to an arbitrary number of blocks (hence allowing $N > 2$). Motivated by the significant range of mathematical models that require the solution of multiple saddle-point systems, we prove that all complex eigenvalues with non-trivial imaginary part (referred to as complex eigenvalues in the following) of the preconditioned system are strictly contained within a circle in the complex plane, with center $1$ and radius bounded by $1$. Furthermore, we develop a sequence of polynomials defined by a particular recurrence relation, the extremal roots of which may be used to describe the end-points of the intervals containing the real (positive) eigenvalues of the preconditioned system. These polynomials are constructed using a related idea to that described in \cite{bergamaschi2026eigenvalue} for the simpler block-diagonal preconditioner. 

{As is well known, eigenvalue information alone cannot completely describe the convergence behavior of Krylov subspace methods such as GMRES. The classical upper bound for the residual norm of GMRES involves the condition number of the eigenvector matrix {in the case of a diagonalizable matrix \cite[Thm. 7.1]{ESW2014}}. Although we do not have theoretical estimates for the condition number of this matrix, we have seen experimentally that its value was frequently modest for the class of problems examined here. Furthermore, numerical results for synthetic test problems in Section \ref{sec:NumExpts} show, unsurprisingly, a correlation between the number of GMRES iterations and the distribution of the real eigenvalues of the preconditioned matrix.
}

The paper is structured as follows: In Section \ref{sec:Preconditioner} we describe the preconditioner and write the eigenvalue problem for the preconditioned matrix. In Section \ref{sec:complex} we characterize the  complex eigenvalues
of the preconditioned matrix. Section \ref{sec:bounds} is devoted to linking the real eigenvalues with the roots of suitable parametric polynomials, while in Section \ref{sec:roots} we describe how to bound such roots.
Section \ref{sec:NumExpts} is devoted to comparing the bounds with the computed eigenvalues on a range of synthetic experiments, with the number of blocks ranging from 3 to 5. The bounds are also verified in a set of
realistic problems arising from PDE-constrained optimization. Some concluding remarks are provided in Section \ref{sec:conc}.

\section{Block-Triangular Preconditioner}
\label{sec:Preconditioner}

We consider the solution the linear system \eqref{Eq1FULL} 
with $\fx$, $\fb$ vectors of consistent length, by an iterative method, preconditioned with the following (inexact) triangular preconditioner:
\begin{equation*}
	\mathcal{P} = \begin{bmatrix}
		\widehat{S}_0 & B_1^\top & & & \\
		& -\widehat{S}_1 & B_2^\top & & \\
		& & \widehat{S}_2 & \ddots &\\
		& & & \ddots & B_N^\top \\
		& & & & (-1)^{N}\widehat{S}_N
\end{bmatrix},
\end{equation*}
where
\begin{equation*}
\begin{array}{ll}
	\widetilde{S}_0 = A_0, \qquad & \widehat{S}_0 \approx \widetilde{S}_0, \\
	\widetilde{S}_k = A_k + B_k\widehat{S}^{-1}_{k-1}B_k^{\top}, \qquad & \widehat{S}_k \approx \widetilde{S}_k,
\end{array}
\end{equation*}
with $\widehat{S}_k$ symmetric positive definite approximations of each $\widetilde{S}_k$. Let
${\mathcal{D}}=\text{diag}\left(\widehat S_0, \widehat S_1, \ldots, \widehat S_N\right).
$
Then, finding the eigenvalues of $\mathcal{P}^{-1}\mathcal{A}$ is equivalent to solving
\begin{equation*} 
\mathcal{D}^{-1/2}\mathcal{A}\mathcal{D}^{-1/2}\fu=\lambda \mathcal{D}^{-1/2}\mathcal{P}\mathcal{D}^{-1/2}\fu, \end{equation*}
or
	\begin{equation} 
		\label{eigMatrix}
	\begin{bmatrix}
	E_0 & R_1^\top &   &        &    \\
	R_1 & -E_1 & R_2^\top &        &        \\
	  &  R_2 &  E_2 &          \ddots &\\
		&      & \ddots          &  \ddots       & R_N^\top \\
		&        &   & R_N &   (-1)^N E_N
\end{bmatrix} \fu =  \lambda 
	\begin{bmatrix}
	I & R_1^\top &   &        &    \\
	   & -I & R_2^\top &        &        \\
	  &    &  I &          \ddots &\\
		&      &          &  \ddots       & R_N^\top \\
		&        &   &   &   (-1)^{N}I 
	\end{bmatrix} \fu,  \qquad \fu = 
\begin{bmatrix} \fu_1\\ \fu_2 \\ \fu_3 \\ \vdots \\ \fu_{N+1} \end{bmatrix}\\\end{equation}
	where $E_k = \widehat S_k^{-1/2} A_k \widehat S_k^{-1/2}$ and $R_k = \widehat S_{k}^{-1/2} B_k \widehat S_{k-1}^{-1/2}$.

		Componentwise, we write the generalized eigenvalue problem \eqref{eigMatrix}
        as
\begin{equation}
\label{dpsp}
\begin{array}{ccccccccc}
	    (E_0-\lambda I)\fu_1 & + & (1-\lambda) R_1^{\top}\fu_2 & & & & & = & \fzero, \\
	    R_1\fu_1 & + &(\lambda I - E_1) \fu_2 & +& (1-\lambda) R_2^{\top}\fu_3 & & & = & \fzero, \\
	    & & R_2\fu_2 & + & (E_2-\lambda I)\fu_3 & + & (1-\lambda) R_3^{\top} \fu_4 & = & \fzero, \\
         & & & \vdots & & \vdots & & \vdots &  \\
	    & & R_{N-1}\fu_{N-1}& + & (-1)^{N-1}(E_{N-1}-\lambda I)\fu_{N} & + & (1-\lambda) R_N^{\top}\fu_{N+1} & = & \fzero, \\
	    & & & & R_N \fu_N & + &(-1)^{N}(E_{N}-\lambda I)\fu_{N+1} & = & \fzero.
    \end{array}
\end{equation}

	\section{Complex eigenvalues}
\label{sec:complex}
    We now wish to characterize the eigenvalues of \eqref{eigMatrix} with non-zero imaginary part. We set $\gamma_E^{(j)} = \dfrac{{\fu_{j+1}^{*}} E_j {\fu_{j+1}}}{\fu_{j+1}^{*} {\fu_{j+1}}} \in [\alpha_E^{(j)}, \beta_E^{(j)}]$, from which it follows that $\alpha_E^{(0)} > 0$ and $\alpha_E^{(j)} \ge 0$ for $j > 0$, and we also define $\varphi_j = \fu_j^* R_j^{\top} \fu_{j+1}$. We then notice that $\fu_1$ must be non-zero, as taking $\fu_1 = \fzero$ and applying forward substitution within the equations of \eqref{dpsp} would imply either that $\lambda = 1$, or that $\fu_2 = \fu_3 = \ldots = \fu_{N+1} = \fzero$, contrary to the definition of an eigenvector. Further, $\fu_2 \ne \fzero$ must hold, otherwise $\lambda \in [\alpha_E^{(0)}, \beta_E^{(0)}]$ is real, or $\fu_1 = \fzero$ which retrieves the previous case.	
    To proceed further, we first assume that all other sub-vectors $\fu_j$, $j \ge 3$, are non-zero. As we will see later, this assumption can be removed while obtaining similar results (see Remark \ref{Rem:1}). 
    Then, premultiplying the $k$-th equation of \eqref{dpsp} by $\fu_k^*$, we obtain
	\begin{equation}
		\label{componentwise}
		\begin{array}{lclclcl}
			&-&(\lambda -\gamma_E^{(0)}) \|\fu_1\|^2  &+&(1-\lambda) \varphi_1& =& \fzero, \\
			\overline \varphi_1       &+& (\lambda  - \gamma_E^{(1)}) \|\fu_2\|^2          &+&(1-\lambda)\varphi_2 &=& \fzero, \\
			 \overline \varphi_2 &-&(\lambda - \gamma_E^{(2)})  \|\fu_3\|^2 &+& (1-\lambda)\varphi_3  &=& \fzero, \\
			& \vdots && \vdots && \vdots & \\
			\overline \varphi_{N-2}  &  +& (-1)^{N-1}(\lambda - \gamma_E^{(N-2)}) \|\fu_{N-1}\| &+& (1-\lambda)\varphi_{N-1}  &=&  \fzero, \\
			\overline \varphi_{N-1}  &  +& (-1)^N (\lambda - \gamma_E^{(N-1)}) \|\fu_N\|^2 &+& (1-\lambda) \varphi_{N}  &=&  \fzero, \\
			 \overline \varphi_N    &+& (-1)^{N+1} (\lambda - \gamma_E^{(N)}) \|{\fu}_{N+1}\|^2 &&& =& \fzero.
		\end{array}
	\end{equation}
	We assume, without loss of generality, that $\|\fu\|^2 = \displaystyle \sum_{k=1}^{N+1} \|\fu_k\|^2 = 1$.
	We rewrite system \eqref{componentwise}, after conjugation of  equations $2, 4, \ldots, 2 \lfloor (N+1)/2 \rfloor$,
	in compact form as \begin{equation}
		\label{sys}
	\mathbf{F}(\fU,\boldsymbol\Phi,\overline{\boldsymbol\Phi}, \lambda) = \fzero, \end{equation}
	where if $N$ is even\footnote{If $N$ is odd the last two components change to
			\[ \begin{bmatrix}
				\overline \varphi_{N-1} - (\lambda  - \gamma_E^{(N-1)}) \|\fu_N\|^2 + (1-\lambda) \varphi_{N}  \\
			\varphi_N     + (\overline \lambda - \gamma_E^{(N)}) \|\fu_{N+1}\|^2  \end{bmatrix}. \]}
	\[ \mathbf{F}(\fU, \boldsymbol\Phi, \overline{\boldsymbol\Phi}, \lambda) = 	\begin{bmatrix}
		- (\lambda - \gamma_E^{(0)}) \|\fu_1\|^2  +(1-\lambda) \varphi_1 \\
			\varphi_1       + (\overline \lambda - \gamma_E^{(1)}) \|\fu_2\|^2          +(1-\overline \lambda)\overline \varphi_2  \\
			 \overline \varphi_2 -(\lambda - \gamma_E^{(2)}) \|\fu_3\|^2 + (1-\lambda)\varphi_3   \\
			\vdots  \\
			\varphi_{N-1}    + (\overline \lambda - \gamma_E^{(N-1)}) \|\fu_N\|^2 + (1-\overline \lambda) \overline \varphi_{N}  \\
			 \overline \varphi_N     - (\lambda  - \gamma_E^{(N)})\|{\fu}_{N+1}  \|^2
		\end{bmatrix}, \quad
	\fU = \begin{bmatrix} \|\fu_1\|^2 \\ \|\fu_2\|^2 \\ \vdots \\\|\fu_{N+1}\|^2\end{bmatrix}, \quad \boldsymbol\Phi=
		\begin{bmatrix} \varphi_1 \\ \varphi_2 \\ \vdots \\\varphi_{N} \end{bmatrix}.\]
The following lemma (implicitly) characterizes the complex eigenvalues:
	\begin{Lemma}
		\label{LemmaComplex}
		Let $\fu_j \ne \fzero$, $j = 1, \ldots, N+1$. The (possibly complex) eigenvalues of 
		\eqref{dpsp} satisfy
		\begin{equation} \label{PolComplex}
		\sum_{j = 1}^{\lfloor (N+1)/2 \rfloor} (1-\lambda) (\overline \lambda-\gamma_E^{(2j-1)})  |\lambda-1|^{2(j-1)} \|\fu_{2j}\|^2 +
			\displaystyle \sum_{j = 0}^{\lfloor N/2 \rfloor}  (\lambda -\gamma_E^{(2j)})  |\lambda-1|^{2j} \|\fu_{2j+1}\|^2 = 0.
			\end {equation}
	\end{Lemma}
	\begin{proof}
		Let us define the diagonal matrix $D(\lambda) = \{d_{i,i}\}$ where the diagonal elements $d_{i,i}$ are defined recursively as
		\begin{align*}
			d_{1,1} & = 1, & & \\
			d_{j,j} &= -(1-\lambda) d_{j-1,j-1},       & \hspace{-10em} 2\le j\le N+1, \quad & j \text{ even}, \\
			d_{j,j} &= -(1-\overline \lambda) d_{j-1,j-1}, & \hspace{-10em} 3\le j\le N+1, \quad & j \text{ odd}.
		\end{align*}
		This provides also a compact definition as 
		\begin{equation*} 
			d_{j,j}  = \begin{cases} -(1-\lambda)^{j/2} (1-\overline \lambda)^{j/2-1}=-(1-\lambda) |\lambda-1|^{j-2}, & j \text{ even}, \\
				(1-\lambda)^{(j-1)/2} (1-\overline \lambda)^{(j-1)/2}=|\lambda-1|^{j-1}, & j \text{ odd}.
			\end{cases}
		\end{equation*}
		System \eqref{sys} is equivalent to $D(\lambda) \mathbf{F}(\fU, \boldsymbol\Phi, \overline{\boldsymbol\Phi}, \lambda) = \fzero$ {if $\lambda \ne 1$}. Componentwise,
		\[ \widehat{\mathbf F} (\fU,\boldsymbol\Phi,\overline{\boldsymbol\Phi}, \lambda) \equiv D(\lambda) \mathbf{F}(\fU,\boldsymbol\Phi,\overline{\boldsymbol\Phi},\lambda) = 	\begin{bmatrix}
			-(\lambda-\gamma_E^{(0)}) \|\fu_1\|^2  +(1-\lambda) \varphi_1 \\
			-(1-\lambda) \varphi_1       - (1-\lambda) (\overline \lambda - \gamma_E^{(1)})\|\fu_2\|^2          -|\lambda-1|^2\overline \varphi_2  \\
			 |\lambda-1|^2 \overline \varphi_2 -|\lambda-1|^2(\lambda - \gamma_E^{(2)})\|\fu_3\|^2 + (1-\lambda)|\lambda-1|^2 \varphi_3   \\
			\vdots  \\
		\end{bmatrix}. \]
			After this scaling, the coefficients of $\varphi_{j}$ (for odd $j$) in two consecutive components of $\widehat{\mathbf{F}}(\fU, \boldsymbol\Phi, \overline{\boldsymbol\Phi}, \lambda)$
			have the same expression with opposite signs, and this also holds for the coefficients of $\overline \varphi_{j}$ (for even $j$). 
			Denoting as $a_j$ the coefficients multiplying $\|\fu_j\|^2$, we easily obtain
		\begin{align*}	
			a_1 &= \gamma_E^{(0)}-\lambda, \\
			a_j &= d_{j,j} (\overline \lambda - \gamma_E^{(j-1)}),\phantom{-} \quad j \text{ even}, \\
			a_j &= -d_{j,j} (\lambda - \gamma_E^{(j-1)}), \quad j \text{ odd}.
		\end{align*}	
	Setting $\fe = \left[ 1, 1, \dots, 1 \right]^{\top}$, the equation
		\[ -\fe^{\top} D(\lambda) \mathbf{F}(\fU,\boldsymbol\Phi,\overline{\boldsymbol\Phi}, \lambda) = \fzero\]
		yields \eqref{PolComplex}.

	\end{proof}
	\begin{Theorem}
    \label{Thm1}
    Let $\fu_j \ne \fzero$, $j = {3}, \ldots, N+1$. If $\alpha_E^{(0)} \|\fu_1\|^2 + \alpha_E^{(1)} \|\fu_2\|^2\ge 1$, no complex eigenvalue can occur.
		Otherwise, all complex eigenvalues $\lambda = \lambda_R  + \emph{\textbf{i}} \, \lambda_I$ satisfy 
		\[ |\lambda -1| < \sqrt{1 - \alpha_E^{(0)} \|\fu_1\|^2- \alpha_E^{(1)} \|\fu_2\|^2}, \qquad \lambda_R >  
		                \frac{\alpha_E^{(0)}\|\fu_1\|^2 + \alpha_E^{(1)} \|\fu_2\|^2}{2}.\]
	\end{Theorem}
	\begin{proof}
                                We will make use of  
the identity $\overline \lambda - |\lambda|^2 = 1- \lambda - |\lambda-1|^2$. We rewrite the first term in  \eqref{PolComplex}
		as
		\begin{align*}
		& \sum_{j = 1}^{\lfloor (N+1)/2 \rfloor}(1-\lambda) (\overline \lambda - \gamma_E^{(2j-1)}) |\lambda-1|^{2(j-1)} \|\fu_{2j}\|^2 
			= \sum_{j = 1}^{\lfloor (N+1)/2 \rfloor}   \left( \overline \lambda - |\lambda|^2 + \gamma_E^{(2j-1)} (\lambda-1)\right)   |\lambda-1|^{2(j-1)}\|\fu_{2j}\|^2\\
			={}&   \sum_{j = 1}^{\lfloor (N+1)/2 \rfloor}  \left(1- |\lambda-1|^2 - \gamma_E^{(2j-1)}  - \lambda(1 - \gamma_E^{(2j-1)} )\right)   |\lambda-1|^{2(j-1)}\|\fu_{2j}\|^2.
\end{align*}
		Then \eqref{PolComplex} may be rewritten as
\begin{align*}
0={}& \underbrace{-\sum_{j = 0}^{\lfloor N/2 \rfloor} \gamma_E^{(2j)} |\lambda-1|^{2j} \|\fu_{2j+1}\|^2+
	  \sum_{j = 1}^{\lfloor (N+1)/2 \rfloor} \left(1- |\lambda-1|^2- \gamma_E^{(2j-1)}\right) |\lambda-1|^{2(j-1)} \|\fu_{2j}\|^2}_
	{\text{\normalsize $\Psi_1$}} \\
&+ \lambda \underbrace{\left(\sum_{j = 0}^{\lfloor N/2 \rfloor}  |\lambda-1|^{2j} \|\fu_{2j+1}\|^2  - \sum_{j = 1}^{\lfloor (N+1)/2 \rfloor}(1-\gamma_E^{(2j-1)})  |\lambda-1|^{2(j-1)}\|\fu_{2j}\|^2\right)}_{\text{\normalsize $\Psi_2$}} \\
	\nonumber \equiv{}& \Psi_1 + (\lambda_R   + \textbf{i}\, \lambda_I)  \Psi_2 = (\Psi_1 + \lambda_R \Psi_2) + \textbf{i} \, \lambda_I  \Psi_2,
\end{align*}	
with $\Psi_1, \Psi_2, \lambda_R, \lambda_I \in \mathbb{R}$. For $\lambda$ to be complex, we must have $\Psi_2 = 0$, which implies also that $\Psi_1 = 0$. In particular,
\begin{equation}
	\label{ReGen}	  0 = \Psi_1 = 
	-\sum_{j = 0}^{\lfloor N/2 \rfloor} \gamma_E^{(2j)} |\lambda-1|^{2j} \|\fu_{2j+1}\|^2+
	  \sum_{j = 1}^{\lfloor (N+1)/2 \rfloor} \left(1- |\lambda-1|^2- \gamma_E^{(2j-1)}\right) |\lambda-1|^{2(j-1)} \|\fu_{2j}\|^2.
\end{equation}	
		Setting $s = |\lambda -1 |^2$, we can view \eqref{ReGen} as a real polynomial equation
			  $Q(s) = 0$, where
			  \begin{align*} Q(s) &=  
				  -\sum_{j = 0}^{\lfloor N/2 \rfloor}  \gamma_E^{(2j)} \|\fu_{2j+1}\|^2 s^{j} 
				  +\sum_{j = 1}^{\lfloor (N+1)/2 \rfloor} \|\fu_{2j}\|^2 s^{j-1} 
			   - s \sum_{j = 1}^{\lfloor (N+1)/2 \rfloor} \|\fu_{2j}\|^2 s^{j-1} -
				  \sum_{j = 1}^{\lfloor (N+1)/2 \rfloor} \gamma_E^{(2j-1)} \|\fu_{2j}\|^2 s^{j-1} \\
				  &= -R_1(s) + P(s) - sP(s) - R_2(s)
				  \end{align*}
        and
		\begin{equation*}
			P(s) = \sum_{j = 1}^{\lfloor (N+1)/2 \rfloor} \|\fu_{2j}\|^2 s^{j-1}, \qquad  
			R_1(s) = \sum_{j = 0}^{\lfloor N/2 \rfloor} \gamma_E^{(2j)} \|\fu_{2j+1}\|^2 s^j, \qquad 
			R_2(s) = \sum_{j = 1}^{\lfloor (N+1)/2 \rfloor} \gamma_E^{(2j-1)} \|\fu_{2j}\|^2 s^{j-1}.
        \end{equation*}
		Notice that if we consider complex eigenvalues we may write $\|\fu_1\|^2, \|\fu_2\|^2 > 0$, as if this does not hold we are reduced to a real eigenvalue case. Given this condition, we have that $P(s) > 0$, $R_1(s) > 0$, and $R_2(s) \ge 0$, $\forall s \ge 0$.
		If $\xi > 0$ is a root of $Q$, it holds that
		\[ \xi = 1 - \frac{R_1(\xi) + R_2(\xi)}{P(\xi)} < 1.\]
          Moreover, $P(\xi) < P(1) = \displaystyle \sum_{j = 1}^{{\lfloor (N+1)/2 \rfloor}} \|\fu_{2j}\|^2 < 1$, 
                and therefore
                \begin{align} \label{ComplexRoots}|\lambda-1|^2 &= \xi 
                = 1 - \frac{R_1(\xi) + R_2(\xi)}{P(\xi)} \\
\nonumber &= 1 -  \frac{\gamma_E^{(0)} \|\fu_1\|^2 + \gamma_E^{(1)} \|\fu_2\|^2  + \sum_{j = 1}^{\lfloor N/2 \rfloor} \gamma_E^{(2j)} \|\fu_{2j+1}\|^2 s^j + \sum_{j = 2}^{\lfloor (N+1)/2 \rfloor} \gamma_E^{(2j-1)} \|\fu_{2j}\|^2 s^{j-1}} {\sum_{j = 1}^{\lfloor (N+1)/2 \rfloor} \|\fu_{2j}\|^2 s^{j-1}} \\ 
\nonumber &< 1 - \gamma_E^{(0)} \|\fu_1\|^2  - \gamma_E^{(1)} \|\fu_2\|^2 \le 1 - \alpha_E^{(0)} \|\fu_1\|^2 - \alpha_E^{(1)} \|\fu_2\|^2.\end{align}
		Finally,
		\[ \lambda_R = \frac{|\lambda|^2 - |\lambda-1|^2 +1}{2} > \frac{1-|\lambda-1|^2}{2} > 
		\frac{\alpha_E^{(0)}\|\fu_1\|^2 + \alpha_E^{(1)} \|\fu_2\|^2}{2}.\]
	\end{proof}

\begin{Remark} \label{Rem:1}
The above results can be extended to the setting where the assumption $\fu_j \ne \fzero$, $j \ge 3$, does not hold. Specifically, Lemma \ref{LemmaComplex} may be modified by defining \[j_{\min} = \min \{j, \ \text{such that} \ \fu_j = \fzero\},\] at which point the eigenvalues of \eqref{eigMatrix} satisfy \eqref{PolComplex} with $j_{\min} -1$ in place of $N+1$. 
Consequently, Theorem \ref{Thm1} still holds true.
\end{Remark}

\begin{Remark}
If, in addition to the assumptions of Theorem \ref{Thm1},
we assume that all matrices $E_j$ with $j$ odd are symmetric positive definite, then we can state a bound for the complex eigenvalues, independent of the norm of the eigenvector. In fact, from \eqref{ComplexRoots} we have
\[\frac{R_1(\xi) + R_2(\xi)}{P(\xi)} 
> \frac{R_2(\xi)}{P(\xi)} 
= \frac{\sum_{j = 1}^{\lfloor (N+1)/2 \rfloor} \gamma_E^{(2j-1)} \|\fu_{2j}\|^2 s^{j-1}} {\sum_{j = 1}^{\lfloor (N+1)/2 \rfloor} \|\fu_{2j}\|^2 s^{j-1}} 
\ge \min\{ \alpha_E^{(2j-1)}, \ j = 1, \ldots, \lfloor (N+1)/2 \rfloor\},\]
which yields that 
\[ |\lambda-1|^2 < 1 - \min\{ \alpha_E^{(2j-1)}, \ j = 1, \ldots, \lfloor (N+1)/2 \rfloor\}, \qquad \lambda_R > \frac 12\min\{ \alpha_E^{(2j-1)}, \ j = 1, \ldots, \lfloor (N+1)/2 \rfloor\}.\]
\end{Remark}
\begin{Remark}
The statement $\Psi_2 = 0$, for a generalized saddle-point system ($N = 1$), simplifies to
                \[ \|\fu_{1}\|^2 - \|\fu_{2}\|^2 + \gamma_E^{(1)} \|\fu_{2}\|^2 = 0 \qquad \Longrightarrow \qquad \|\fu_{1}\|^2 = (1 - \gamma_E^{(1)}) \|\fu_{2}\|^2,\]
showing that complex eigenvalues can occur only if $\gamma_E^{(1)} < 1$.
                We then rewrite  \eqref{ComplexRoots} as
                \begin{align*} |\lambda-1|^2 &= 1 - \frac{\gamma_E^{(0)} \|\fu_1\|^2  + \gamma_E^{(1)}\|\fu_2\|^2}{\|\fu_2\|^2} \\
                &=1 - \frac{\gamma_E^{(0)} (1 - \gamma_E^{(1)}) \|\fu_{2}\|^2 + \gamma_E^{(1)}\|\fu_2\|^2}{\|\fu_2\|^2}
                = (1-\gamma_E^{(0)}) (1-\gamma_E^{(1)}) \le
                                (1-\alpha_E^{(0)}) (1-\alpha_E^{(1)}), \end{align*}
since $\gamma^{(0)}_E < 1$ must clearly hold as well, otherwise we would have a contradiction to an absolute value being non-negative, or $\lambda = 1 \in \R$. This bound is slightly tighter than that in \cite[Theorem 2]{Sim04}, and coincides with the bound in \cite[Theorem 5]{BergaNLAA10}.
        \end{Remark}
        \begin{Remark}
        The expression in \eqref{ComplexRoots} agrees with that in \cite[Eq. 2.20]{BerFerMar26}. In fact, for a double saddle-point linear system ($N=2$),
        setting $\rho  = 1-|\lambda-1|^2$, \eqref{ComplexRoots}  simplifies as follows:
        \begin{align*}  \rho = \frac{R_1(\xi) + R_2(\xi)}{P(\xi)} &= \frac{\gamma_E^{(0)} \|\fu_1\|^2 + \gamma_E^{(1)} \|\fu_2\|^2+ \gamma_E^{(2)} |\lambda-1|^2 \|\fu_3\|^2} {\|\fu_2\|^2} \\ 
         &= \frac{\gamma_E^{(0)} \|\fu_1\|^2 + \gamma_E^{(1)} \|\fu_2\|^2 + \gamma_E^{(2)} \|\fu_3\|^2 - \gamma_E^{(2)} \rho \|\fu_3\|^2}{\|\fu_2\|^2}. \end{align*}
         Solving for $\rho$ yields
        \[ \rho = \frac{\gamma_E^{(0)} \|\fu_1\|^2 + \gamma_E^{(1)} \|\fu_2\|^2+ \gamma_E^{(2)} \|\fu_3\|^2}{\|\fu_2\|^2 + \gamma_E^{(2)} \|\fu_3\|^2},\]
which is equivalent to \cite[Eq. 2.20]{BerFerMar26}.
\end{Remark}


	\section{Intervals of real eigenvalues}
\label{sec:bounds}


To analyse the real eigenvalues of the preconditioned system, we now return to \eqref{dpsp}. We note that the matrices $R_kR_k^{\top}$ are all symmetric positive definite. We define two indicators $\gamma_E^{(k)}$ and $\gamma_R^{(k)}$ using the Rayleigh quotient, and denote as $\fw$ a generic real non-zero vector of consistent dimension, to write:
\[
\begin{array}{llll}
	\alpha_E^{(k)} \equiv \lambda_{\min }\left(E_k\right), & \beta_E^{(k)} \equiv \lambda_{\max }\left(E_k\right), & \gamma_E^{(k)}\left(\fw\right)=\dfrac{\fw^{\top} E_k \fw}{\fw^{\top} \fw} \in\left[\alpha_{E}^{(k)}, \beta_{E}^{(k)}\right], 
    & k=0,\dots,N, \\[.6em]
\alpha_R^{(k)} \equiv \lambda_{\min }\left(R_k R_k^{\top}\right), & \beta_R^{(k)} \equiv \lambda_{\max }\left(R_k R_k^{\top}\right), & \gamma_R^{(k)}\left(\fw\right)=\dfrac{\fw^{\top} R_k R_k^{\top} \fw}{\fw^{\top} \fw} \in\left[\alpha_{R}^{(k)}, \beta_{R}^{(k)}\right], 
& k=1,\dots,N. \\
\end{array}
\] 
We define the two vectors of parameters:
\[
\boldsymbol{\gamma}_E=\left[\gamma_E^{(0)},\dots,\gamma_E^{(N)}\right],\quad\quad \boldsymbol{\gamma}_R=\left[\gamma_R^{(1)},\dots,\gamma_R^{(N)}\right], 
\]
and introduce a sequence of parametric polynomials depending on $\boldsymbol{\gamma}_E$, $\boldsymbol{\gamma}_R$ and on $\lambda \in \R$:
\begin{Definition} 
\label{defV}
    \begin{align} \nonumber V_0(\lambda,\boldsymbol{\gamma}_R,\boldsymbol{\gamma}_E) &= 1,\\ 
	\nonumber V_1(\lambda,\boldsymbol{\gamma}_R,\boldsymbol{\gamma}_E) &= \lambda -\gamma_E^{(0)},\\ 
    \label{eq2} V_{k+1}(\lambda,\boldsymbol{\gamma}_R,\boldsymbol{\gamma}_E) &= (\lambda - \gamma_E^{(k)})V_k(\lambda,\boldsymbol{\gamma}_R,\boldsymbol{\gamma}_E)- (\lambda -1)\gamma_R^{(k)}V_{k-1}(\lambda,\boldsymbol{\gamma}_R,\boldsymbol{\gamma}_E), \quad k \ge 1. \end{align}
\end{Definition}
\noindent
Relatedly, we define a sequence of real matrix-valued functions:
\begin{Definition}
    \[
    \begin{aligned}
        Y_1(\lambda) &= \lambda I - E_0, \\
	    Y_{k+1}(\lambda) &= (\lambda-1) R_k Y_{k}(\lambda)^{-1}R_k^{\top} + (-1)^{k}\left(\lambda I - E_k\right), \quad k \ge 1.
    \end{aligned}
    \]
\end{Definition}
We now introduce a technical lemma, the proof of which can be found in \cite[Lem. 2.1]{BMPP_COAP25}.
\begin{Lemma}\label{lemma1}
    \textit{Let $Y$ be a real symmetric matrix-valued function defined in $F\subset\mathbb{R}$ and}
    \[
        0\notin[\min\{\sigma(Y(\zeta))\},\max\{\sigma(Y(\zeta))\}]\quad \quad \textit{for all }\ \zeta\in F.
    \]
    \textit{Then, for arbitrary real $\fs\neq\fzero$, there exists a real vector $\fv\neq\fzero$ such that}
    \[
        \frac{\fs^{\top}Y(\zeta)^{-1}\fs}{\fs^{\top}\fs}= \frac{1}{\gamma_Z} \qquad \textit{with }\ \gamma_Z = \frac{\fv^{\top}Y(\zeta)\fv}{\fv^{\top}\fv}.
    \]
\end{Lemma}

\medskip

\noindent
   \textbf{Notation}.
 We use the notation
$\mathcal{I}_k$ to denote the union of the minimal closed
intervals containing all real roots of polynomials of the form $V_k$ 
over the
valid range of $\boldsymbol\gamma_E$, $\boldsymbol\gamma_R$. 
We also set $\mathcal{I}_{\leq k} := \bigcup_{j = 1}^k \mathcal{I}_k$.
Throughout our derivations, we frequently use $V_k$ as a shorthand for $V_k(\lambda,\boldsymbol{\gamma}_R,\boldsymbol{\gamma}_E)$.

\medskip

The next lemma will be used in the proof of the subsequent Theorem \ref{theorem1}.

\begin{Lemma}\label{lemma2}
    For every real vector $\fu\neq\fzero$ {and $\lambda\in\mathbb{R}$}, there is a choice of $\boldsymbol{\gamma}$ for which 
    \[
	    \frac{\fu^{\top}Y_{k+1}(\lambda)\fu}{\fu^{\top}\fu}=(-1)^{k} \frac{V_{k+1}(\lambda)}{V_k(\lambda)}\quad\quad \textit{for all }\ \lambda\notin  \mathcal{I}_{\leq k}. \]
\end{Lemma}
\begin{proof}
	This is shown by induction. For $k=0$ we have $\displaystyle\frac{\fu^{\top}Y_1(\lambda)\fu}{\fu^{\top}\fu}=\lambda -\gamma_E^{(0)}=\frac{V_1(\lambda)}{V_0(\lambda)}$ for all $\lambda\in\mathbb{R}$. If $k\geq1$ the condition $\lambda\notin{\mathcal{I}}_k$, together with the inductive hypothesis $\displaystyle\frac{\fu^{\top}Y_{k}(\lambda)\fu}{\fu^{\top}\fu}=(-1)^{k-1}\frac{V_{k}(\lambda)}{V_{k-1}(\lambda)}$, implies invertibility of $Y_k(\lambda)$. Moreover, this is equivalent to the condition $0\notin[\min\{\sigma(Y(\zeta))\},\max\{\sigma(Y(\zeta))\}]$ that guarantees the applicability of Lemma \ref{lemma1}. Therefore, we can write 
\begin{align}
	\nonumber \frac{\fu^{\top}Y_{k+1}(\lambda)\fu}{\fu^{\top}\fu}&\underset{\color{white} \scalebox{0.7}{$\fw$} = R^{\top}_k \scalebox{0.7}{$\fu$}}{=}(\lambda-1)\frac{\fu^{\top}R_k Y_{k}(\lambda)^{-1}R_k^{\top}\fu}{\fu^{\top}\fu} + (-1)^k(\lambda-\gamma_E^{(k)})\\
	\label{eq3} &\underset{\scalebox{0.7}{$\fw$} = R^{\top}_k \scalebox{0.7}{$\fu$}}{=}(\lambda-1)\frac{\fw^{\top}Y_{k}(\lambda)^{-1}\fw}{\fw^{\top}\fw}\gamma_R^{(k)} + (-1)^k(\lambda-\gamma_E^{(k)}).
\end{align}
We then apply Lemma \ref{lemma1} and the inductive hypothesis to write 
\[
	\frac{\fw^{\top}Y_{k}(\lambda)^{-1}\fw}{\fw^{\top}\fw}=(-1)^{k-1}\frac{V_{k-1}(\lambda)}{V_{k}(\lambda)}.
\]
Substituting into \eqref{eq3} and using relation \eqref{eq2}, we obtain 
	\begin{align*}(\lambda-1)\frac{\fw^{\top}Y_{k}(\lambda)^{-1}\fw}{\fw^{\top}\fw}\gamma_R^{(k)} + (-1)^{k}(\lambda -\gamma_E^{(k)}) &=
		(-1)^{k-1}(\lambda-1) \gamma_R^{(k)}\frac{V_{k-1}(\lambda)}{V_{k}(\lambda)} + (-1)^{k}(\lambda-\gamma_E^{(k)}) \\
		&= \frac{(-1)^{k} \left(-(\lambda -1) \gamma_R^{(k)} V_{k-1}(\lambda) + (\lambda - \gamma_E^{(k)})V_k(\lambda)\right)}{V_k(\lambda)} \\
		&= (-1)^{k} \frac{V_{k+1}(\lambda)}{V_k(\lambda)}.
	\end{align*}
\end{proof}

\begin{Theorem}\label{theorem1}
    {Any {real eigenvalue} $\lambda$ of $\mathcal{P}^{-1}\mathcal{A}$ is located in $\mathcal{I}_{\leq N+1}$}.
\end{Theorem}

\begin{proof}
The proof is carried out through induction on $k$, that is for every $k\leq N+1$ either
\[
\text{(i)}\ \lambda\in{\mathcal{I}}_k \quad\quad\text{or}\quad\quad\text{(ii)}\ \fu_k={(1-\lambda)}Y_{k}(\lambda)^{-1}R_k^{\top}\fu_{k+1}.
\]
Let $\fu=\left[\fu_1^{\top},\dots,\fu_{N+1}^{\top}\right]^{\top}$ be a real eigenvector of \eqref{eigMatrix}. Assuming that $\lambda\notin [\alpha_E^{(0)}, \beta_E^{(0)}]$, then $Y_1(\lambda)$ is invertible. From the first row of \eqref{dpsp} we obtain 
\begin{equation}
	(\lambda I -E_0)\fu_1 = (1-\lambda) R_1^{\top} \fu_2\quad\Longrightarrow\quad Y_1(\lambda)\fu_1 = (1-\lambda) R_1^{\top} \fu_2,
\label{eq4}
\end{equation}
whereupon inserting \eqref{eq4} into the second row of \eqref{dpsp} yields 
\begin{equation}
	\left((1-\lambda) R_1Y_1(\lambda)^{-1}R_1^{\top}{+}\lambda I {-}E_1\right)\fu_2=(\lambda-1) R_2^{\top}\fu_3 \quad\Longrightarrow\quad Y_2(\lambda)\fu_2={(1-\lambda)} R_2^{\top}\fu_3.
\label{eq5}
\end{equation}
Pre-multiplying the left-hand side of \eqref{eq5} by $\frac{\fu_2^{\top}}{\fu_2^{\top}\fu_2}$, we obtain that
\begin{equation*}
\begin{aligned}
	-\frac{\fu_2^{\top}Y_2(\lambda)\fu_2}{\fu_2^{\top}\fu_2} &\underset{\phantom{\scalebox{0.7}{$\fs$}=R_1^{\top} \scalebox{0.7}{$\fu_2$}}}{=}\frac{\fu_2^{\top}((1-\lambda) R_1Y_1(\lambda)^{-1}R_1^{\top})\fu_2}{\fu_2^{\top}\fu_2} + \lambda - \frac{\fu_2^{\top}E_1\fu_2}{\fu_2^{\top}\fu_2}\\
	&\underset{\scalebox{0.7}{$\fs$}=R_1^{\top} \scalebox{0.7}{$\fu_2$}}{=}(1-\lambda)\frac{\fs^{\top}Y_1(\lambda)^{-1}\fs}{\fs^{\top}\fs}\gamma_R^{(1)} + \lambda - \gamma_E^{(1)}.
\end{aligned}
\end{equation*}
By virtue of \eqref{eq2}, this expression equals
\[
	(1-\lambda)    \frac{V_0(\lambda)}{V_1(\lambda)}\gamma_R^{(1)} + \lambda - \gamma_E^{(1)}= \frac{V_2(\lambda)}{V_1(\lambda)}.
\]
If $\fu_3=\fzero$, then $\lambda \in \R$ being an eigenvalue implies that $V_2(\lambda,\boldsymbol{\gamma}_E, \boldsymbol{\gamma}_R)=0$. Otherwise, if $\lambda\notin{\mathcal{I}}_2$ then $Y_2(\lambda)$ is invertible.
Assume now the inductive hypothesis holds for $k-1$. If $\lambda\notin{\mathcal{I}}_{k-1}$, then $Y_{k-1}(\lambda)$ is definite and invertible. We can write
\[
	\left((\lambda-1) R_{k-1}Y_{k-1}(\lambda)^{-1}R_{k-1}^{\top}+(-1)^{k-1}(\lambda I-E_{k-1}) 
    \right)
    \fu_k=(1-\lambda) R_k^{\top}\fu_{k+1}.
\]
Now, if $\lambda \in \R$ is an eigenvalue with $\fu_{k+1}=\fzero$, from
\[
0=\frac{\fu_k^{\top}Y_{k}(\lambda)\fu_k}{\fu_k^{\top}\fu_k}=(-1)^{k-1}\frac{V_k(\lambda)}{V_{k-1}(\lambda)}
\]
we have that $V_{k}(\lambda,\boldsymbol{\gamma}_E, \boldsymbol{\gamma}_R)=0$ and hence $\lambda\in{\mathcal{I}}_k$. Otherwise, if $\lambda\notin{\mathcal{I}}_k$ we may write
\[
	\fu_k=(1-\lambda)Y_{k}(\lambda)^{-1} R_k^{\top}\fu_{k+1}.
\]
The induction process ends for $k=N+1$. In this case we have that 
\[
Y_{N+1}\fu_{N+1}=\fzero,
\]
and the condition $\lambda\in{\mathcal{I}}_{N+1}$ must be achieved, noticing that $\fu_{N+1}=\fzero$ would imply that also $\fu_N=\ldots=\fu_1=\fzero$ contradicting the definition of an eigenvector.
\end{proof}

\section{Bounds for the extremal real roots of $V_{k+1}$}
\label{sec:roots}
In the previous section we have seen that there is a strict relation between eigenvalues of the preconditioned matrix and
the roots of the sequence of polynomials \eqref{eq2}.
In this section we will characterize the real roots. 
We start by showing that the real roots of $V_k(\lambda)$ (and hence the real eigenvalues of \eqref{dpsp}) are strictly positive. In fact we show by induction that the sign of $V_k(z)$ is $(-1)^k$, for
all $z \le 0$. We start by observing that for $z \le 0$, $V_0(z) = 1 > 0$ and $V_1(z) = z-\gamma_E^{(0)}  < 0$. Assuming that $V_{2j}(z) > 0$, $V_{2j+1}(z) < 0$ for $j \ge 0$, we have 
\begin{align*}
V_{2j+2}(z) &= (z-\gamma_E^{(2j+1)}) V_{2j+1}(z) - (z-1)\gamma_R^{(2j+1)}V_{2j}(z) > 0 , \\ 
V_{2j+3}(z) &= (z-\gamma_E^{(2j+2)}) V_{2j+2}(z) - (z-1)\gamma_R^{(2j+2)} V_{2j+1}(z) < 0.
\end{align*}

Now the question is: for which values of the parameters $\gamma^{(i)}_E$, $\gamma^{(i)}_R$ are the extremal values of the real roots attained?
We begin by showing that the extremal real roots of $V_{k+1}$ are attained when the parameters take extremal values.
       
Let $\xi_1^{(k+1)}, \ldots, \xi_t^{(k+1)}$, $t \leq k+1$, 
denote the $t$ real roots of the polynomial $V_{k+1}$ for a particular combination of the parameters. 
Taking separately one of the parameters, $\gamma_*$ we can write $\gamma_*$ as a convex combination of its extremal values, that is $\gamma_* = \alpha \gamma_*^{\min} + (1-\alpha) \gamma_*^{\max}$ for $\alpha \in [0,1]$, and write, for two polynomials $s_1$ and $s_2$,
\begin{align*}
    V_{k+1} (\lambda, \gamma_*) &= s_1(\lambda)  + \gamma_* s_2(\lambda) =  s_1(\lambda)  + (\alpha \gamma_*^{\min} + (1-\alpha) \gamma_*^{\max})  s_2(\lambda)  \\
    &= \alpha(s_1(\lambda)  + \gamma_*^{\min} s_2(\lambda)) + (1-\alpha) (s_1(\lambda) + \gamma_*^{\max}  s_2(\lambda))
    = \alpha V_{k+1}(\lambda, \gamma_*^{\min}) + (1-\alpha) V_{k+1}(\lambda, \gamma_*^{\max}).
\end{align*}
Therefore,
\[ 0 = V_{k+1} (\xi_j^{(k+1)}, \gamma_*) =   \alpha V_{k+1}(\xi_j^{(k+1)}, \gamma_*^{\min}) + (1-\alpha) V_{k+1}(\xi_j^{(k+1)}, \gamma_*^{\max}). \]
Hence if $\alpha \notin \{0,1\}$, that is $\gamma_* \notin \{\gamma_*^{\min},\gamma_*^{\max}\}$, the values of $\eta_1 = V_{k+1}(\xi_j^{(k+1)}, \gamma_*^{\min})$ and $\eta_2 = V_{k+1}(\xi_j^{(k+1)}, \gamma_*^{\max})$ either take opposite signs, or both are zero. 
In the first case, for the minimum real root $\xi_1^{(k+1)}$, we select $\gamma \in \{\gamma_*^{\min}, \gamma_*^{\max}\}$ corresponding to the $\eta_* \in \{\eta_1, \eta_2\}$ with sign \emph{opposite} to that of $V_{k+1}(0)$, which is $(-1)^{k+1}$, to obtain the lower bound for this root, depending on $\gamma_*$. Regarding the maximum real root $\xi_t^{(k+1)}$, we choose $\gamma \in \{\gamma_*^{\min}, \gamma_*^{\max}\}$ corresponding to the $\eta_* \in \{\eta_1, \eta_2\}$ with \emph{negative} value (the opposite of the sign $V_{k+1}$ has at infinity) to obtain an upper bound for this root. 
If $\eta_1 = \eta_2 = 0$, then $V_{k+1} (\xi_j^{(k+1)}, \gamma_*)$ as a function of $\gamma_*$ is a linear function with two different roots so it must be identically zero. In particular we can take either of the extremal values of the parameter $\gamma_*$.
\begin{Lemma}
\label{increasing}
Let $\xi_{\emph{right}}^{(k+1)}$ be the maximum real (i.e., rightmost) root of $V_{k+1}(\lambda)$, and let $j \leq k$ be the smallest index such that $\xi_\emph{right}^{(j)} > 1$, if such an index exists. Then it holds that
\begin{equation*}
1 < \xi_\emph{right}^{(j)} < \xi_\emph{right}^{(j+1)} < \ldots <  \xi_\emph{right}^{(k+1)}.
\end{equation*}
That is, if there exists an index $j$ such that $\xi_\emph{right}^{(j)} > 1$, the value of $\xi_\emph{right}^{(*)}$ increases thereafter, every time the index is increased.
\end{Lemma}
\begin{proof}
By induction. Assume true the inductive hypothesis up to index $k$, in particular that $\xi_\text{right}^{(k)} > \xi_\text{right}^{(k-1)}$. This implies $V_{k-1}(\xi_\text{right}^{(k)}) > 0$,
then
\[ V_{k+1}(\xi_\text{right}^{(k)}) = (1-\xi_\text{right}^{(k)})\gamma_R^{(k)} V_{k-1}(\xi_\text{right}^{(k)}) < 0,\]
from which, recalling that $\displaystyle \lim_{\lambda  \to +\infty} V_{k+1}(\lambda) = +\infty$, we have $\xi_\text{right}^{(k+1)} > \xi_\text{right}^{(k)}$.
\end{proof}

\begin{Corollary}
  Using the notation of Lemma \ref{increasing}, if  $\xi^{(k+1)}_{\textnormal{right}} > 1$ then $
  \xi^{(k+1)}_{\textnormal{right}} > \xi^{(l)}_{\textnormal{right}}$ for all $1 \leq l \leq k$.
\end{Corollary}
\begin{proof}
  Let $j$ be the smallest index such that  $\xi^{(j)}_{\text{right}} > 1$. 
  Then, if $l<j$, $\xi^{(l)}_{\text{right}} \leq 1$ and $\xi^{(k+1)}_{\text{right}} > 1$, so 
  $\xi^{(k+1)}_{\text{right}} > \xi^{(l)}_{\text{right}}$. If $l \ge j$ then the statement is 
true as a consequence of Lemma \ref{increasing}.
\end{proof}




\newcommand{\vV}{\widetilde V}
\subsection{Largest root}
We now consider the problem of finding for which combination of the extremal values of the parameters 
the largest real root of the polynomial $V_{k+1}(\lambda, \bgamma_E, \bgamma_R)$ takes its maximum value.

\begin{Definition}
\label{defW}
We define a new set of polynomials as
\[ W_0(\lambda) = 1, \qquad W_1(\lambda) = \lambda - \gamma_E^{(k)}, \qquad
 W_{j+1} = (\lambda - \gamma_E^{(k-j)})W_{j} - (\lambda - 1)\gamma_R^{(k+1-j)}
  W_{j-1},\quad{1 \le j \le k}.\]
\end{Definition}
From this definition, it is clear that the $W_j$ satisfy 
the same recursion as the $V_j$ except 
that  $\gamma_E^{(j)}$ is replaced by 
$\gamma_E^{(k-j)}$ and $\gamma_R^{(j)}$ is 
replaced by $\gamma_R^{(k+1-j)}$.
One may also easily show by induction that 
$W_j$ depends only on  $\gamma_E^{(l)} $
for $k \ge l \ge k+1-j$ and on $
\gamma_R^{(l)}$ for $k \geq l \geq k+2-j$.


\begin{Lemma}\label{VkWkRecurrence}
  For all $1 \leq j \leq k$, it holds that 
  \begin{equation}
    \label{dec}
    V_{k+1} = V_jW_{k+1-j} - (\lambda - 1) \gamma_R^{(j)}
    V_{j-1}W_{k-j}.
  \end{equation}
\end{Lemma}
\begin{proof}
  We will show this by induction. 
  The base step is given by $j = k$:
   \[
     V_{k+1} = (\lambda - \gamma_E^{(k)})
     V_k - (\lambda - 1) \gamma_R^{(k)} V_{k-1}
    = V_{k}W_1 - (\lambda - 1) \gamma_R^{(k)}
    V_{k-1}W_0
  .\] 
  Now suppose by induction hypothesis 
  that the statement holds for index $j+1 \leq k$.
  Then 
   \begin{align*}
     V_{k+1} &= V_{j+1}W_{k-j} - (\lambda - 1) \gamma_R^{(j+1)} 
     V_{j}W_{k-1-j} \\ &=
     \left[
       (\lambda  - \gamma_E^{(j)})V_j - (\lambda - 1) \gamma_R^{(j)} V_{j-1}
     \right]W_{k-j} - (\lambda - 1) \gamma_R^{(j+1)}
     V_jW_{k-1-j} \\ &= 
     V_j \left[
       (\lambda - \gamma_E^{(j)})W_{k-j} - (\lambda - 1) \gamma_R^{(j+1)}
     W_{k-1-j} 
   \right] - (\lambda - 1) \gamma_R^{(j)}
     V_{j-1}W_{k-j} \\ &= 
     V_jW_{k+1-j} - (\lambda - 1) \gamma_R^{(j)}
     V_{j-1}W_{k-j}.
  \end{align*}
\end{proof}
\begin{Corollary}
  \label{V=W}
    $V_{k+1} = W_{k+1}$.
\end{Corollary}
\begin{proof}
We start from Lemma \ref{VkWkRecurrence} with $j = 1$ and obtain 
\begin{align*}
V_{k+1} &= V_1W_k - (\lambda - 1) \gamma_R^{(1)}V_0W_{k-1} \\ &= 
(\lambda - \gamma_E^{(0)})W_k - (\lambda - 1) \gamma_R^{(1)}W_{k-1} = W_{k+1}.
\end{align*}
\end{proof}

\begin{Remark}
  Corollary \ref{V=W} tells us that given a choice of parameters $\gamma_E^{(0)}, 
  \dots, \gamma_E^{(k)}$ and $\gamma_R^{(1)}, \dots, \gamma_R^{(k)}$, the polynomials 
  $V_{k+1}$ and  $W_{k+1}$ obtained from the recursions in Definitions \ref{defV} and \ref{defW}, respectively, are the same. It does not imply that $V_j = W_j$ for $j < k+1$, since these polynomials do not depend on the same parameters.
\end{Remark}

\vspace{4mm}

We now assume that the largest real root of $V_{k+1}$ is larger than $1$, otherwise we consider $1$ as its upper bound. 
Under this assumption, we will now show that the largest root of $V_{k+1}$ occurs when all parameters satisfy  $\gamma_E^{(j)} = \beta_E^{(j)}$, $\gamma_R^{(j)} = \beta_R^{(j)}$, for $j \le k$.
We denote by $e_j \in \{\alpha_E^{(j)}, \beta_E^{(j)}\}$ one choice of the extremal values of $\gamma_E^{(j)}$ and, analogously,
$r_j \in \{\alpha_R^{(j)}, \beta_R^{(j)}\}$; moreover, 
we denote by $e^{*}_j$  and $r_j^{*}$ the other extremal value of the parameter.
We now consider the polynomial {$V_{k+1}$} with $\gamma_E^{(j)} \equiv e_j$, $\gamma_R^{(j)} \equiv r_j$, and also define {$\vphantom{V}^{E}V_{k+1}^{j}$} for $0 \leq j \leq k$ as the polynomial 
  {$V_{k+1}$} with $\gamma_E^{(j)} \equiv e_j^{*}$  instead of $e_j$.
   Similarly, we define 
   {$\vphantom{V}^{R}V_{k+1}^{j}$} for $1 \leq j \leq k$ as the 
   polynomial {$V_{k+1}$} with $\gamma_R^{(j)}= r_j^{*}$ instead of $r_j$.
   Then it immediately follows from the 
   recursion for $V_k$ in \eqref{eq2} that 
    \begin{align}
	\label{EV}  \vphantom{V}^{E}V^{j}_{j+1} &=
	  V_{j+1} + (e_j - e_j^{*})
	  V_j & \text{for } \, 0 \leq j \leq k, \\ 
      \nonumber \vphantom{V}^{R}V_{j+1}^{j} &=
	  V_{j+1} - (\lambda - 1) (r_j^{*} - r_j) V_{j-1} & \text{for } \, 1 \leq j \leq k. 
   \end{align}
 We have the following relations: 
 \begin{Lemma} 
 \label{Lemma5}
  \begin{align}
\label{EV_1}	\vphantom{V}^{E}V_{k+1}^{j} &= V_{k+1} +  
	(e_j - e_j^{*}) V_j W_{k-j} & 
    \textnormal{for } \, 0 \leq j \leq k, \\ 
\label{RV_1}	\vphantom{V}^{R}V_{k+1}^{j} &= V_{k+1} - (\lambda - 1) (r_j^{*} - r_j) V_{j-1} W_{k-j}
& \textnormal{for } \, 1 \leq j \leq k.
  \end{align}
\end{Lemma}
\begin{proof}
The case $j = k$ corresponds exactly to the previous formulas, so we can assume $j < k$ from now on. To prove \eqref{EV_1} we start from \eqref{dec}, replacing $j$ with  $j+1$:
   \[
     V_{k+1} = V_{j+1}W_{k-j} - (\lambda - 1) r_{j+1}
    V_{j}W_{k-1-j}.
  \] Then we notice that the only term 
  in the right-hand side that depends on 
  $\gamma_E^{(j)}$ is $V_{j+1}$. Thus,
  \begin{align*}
    \vphantom{V}^{E}V^{j}_{k+1} &= 
    \vphantom{V}^{E}V^{j}_{j+1}W_{k-j} - (\lambda - 1) r_{j+1}
    V_{j}W_{k-1-j} \\ &= 
    V_{k+1} + (
      \vphantom{V}^{E}V_{j+1}^{j} - 
      V_{j+1}
  )W_{k-j} \\ &= 
  V_{k+1} + (e_j - e_j^{*})V_jW_{k-j}.
  \end{align*}
To show \eqref{RV_1}, we start from 
  \eqref{dec}:
\[
    V_{k+1} = V_jW_{k+1-j} - (\lambda - 1) r_j 
    V_{j-1}W_{k-j}
\] and notice that, apart from the constant 
$r_j$, there is no other term on the 
right-hand side that depends on $r_j$.
Thus, 
 \begin{align*}
   \vphantom{V}^{R}V_{k+1}^{j} &= 
    V_jW_{k+1-j} - (\lambda - 1) r_j^{*} 
   V_{j-1}W_{k-j} \\ &= 
   V_{k+1} - (\lambda - 1) (r_j^{*} - r_j)
   V_{j-1}W_{k-j}.
\end{align*}
\end{proof}
\begin{Theorem}
Let $e_j$, $r_j$ be parameters that give the maximum possible real root $\xi$ of $V_{k+1}$ and suppose $\xi > 1$. 
This is achieved by letting all parameters assume the maximum value.
\end{Theorem}
\begin{proof}
    Let $\xi$ be the largest real root of $V_{k+1}$ and let $e_j$, $r_j$ be the choice of parameters that realize it. We start from \eqref{EV_1} and evaluate it at $\xi$, obtaining
    \[ \vphantom{V}^{E}V_{k+1}^{j}(\xi) = (e_j - e_j^*)V_j(\xi) W_{k-j}(\xi).\]
We know from Lemma \ref{increasing} that $\xi$ is larger than all the real roots of $V_j$ and so 
$V_j(\xi) > 0$. 
From Lemma \ref{increasing} and the fact that the $
V_j$ and $W_j$ satisfy the same recursion (except with different parameters), it follows that 
the largest roots of the $W_j$ are also increasing. From Corollary \ref{V=W} we know that $\xi$ 
is the largest root of $W_{k+1}$. From this we see that $W_{k-j}(\xi) > 0$. 
Finally, we are assuming that $\xi$ is the largest possible root (no other choices of parameters can give a larger root). In particular 
$\vphantom{V}^{E}V_{k+1}^{j}(\xi) \ge 0$, with equality holding only if $\xi$ does not depend on $\gamma_E^{(j)}$. Thus we have $e_j - e_j^* \geq 0$, which shows that
$e_j = \beta_E^{(j)}$. We reason similarly for the parameters $\gamma_R^{(j)}$. We start from \eqref{RV_1} and evaluate it at $\xi$ obtaining
\[ \vphantom{V}^{R}V_{k+1}^{j}(\xi) = (1-\xi) (r_j^* - r_j) V_{j-1}(\xi) W_{k-j}(\xi).\]
The same reasoning as above shows that $W_{k-j}(\xi)$ and $V_{j-1}(\xi)$ are positive, and $\vphantom{V}^{R}V_{k+1}^{j}(\xi)$ is non-negative. Furthermore, $1-\xi$ is negative since $\xi > 1$. But then $r_j^* - r_j \leq 0$, which shows that $r_j = \beta_R^{(j)}$.
\end{proof}

\subsection{Smallest root}
To provide a lower bound for the smallest real root of $V_{k+1}$ we assume that this is smaller than $1$. Otherwise we consider $1$ as the lower bound.

\subsubsection{Choice of $\gamma_E^{(j)}$}

\begin{Lemma}
  \label{Wn0}
  Let $\xi$ be a {real} root of  $V_{k+1}$ 
  that lies outside of  $\mathcal{I}_{\leq k}$ 
  and is not equal to $1$. 
  Then  $W_j(\xi) \neq 0$ for all 
  $j \le k$. 
\end{Lemma}
\begin{proof}
Suppose $W_k(\xi) = 0$. Then since $W_{k+1}(\xi) = V_{k+1}(\xi) = 0 $ we see that the value at $\xi$ of two consecutive $W_j$ is zero. Applying the recursion in Definition \ref{defW} (and using the fact that $\xi$ is not equal to $1$) we 
see that $W_j(\xi) = 0$ for all $j \le k+1$. 
This then yields a contradiction since $W_0 = 1$. Now suppose $W_l(\xi) = 0$ for $0 \leq l \leq k-1$. 
We begin by setting $j = k-l$ in \eqref{dec}:
  \[
    V_{k+1} = V_{k-l}W_{l+1} - 
    (\lambda - 1) \gamma_R^{(k-l)} V_{k-1-l}W_{l}
  \] and evaluate at $\xi$ to obtain
  \[
    V_{k-l}(\xi)W_{l+1}(\xi) =
    (\xi - 1) \gamma_R^{(k-l)}
    V_{k-1-l}(\xi)W_{l}(\xi)
  = 0. 
\] 
Since $\xi$ lies outside of $\mathcal{I}_{\leq k}$, it follows that  $V_{k-l}(\xi) \neq 0$ and so $W_{l+1}(\xi) = 0$. Again we have that the value at $\xi$ of two consecutive $W_j$ is zero, which is a contradiction.

\end{proof}
\begin{Theorem} \label{choice}
Let $e_j$, $r_j$ be parameters that give the smallest possible real root $\xi$ of $V_{k+1}$. 
Suppose that $\xi$ lies outside of $\mathcal{I}_{\leq k}$ and $\xi < 1$. 
Then this is achieved by letting the parameters $\gamma_E^{(j)}$ alternate between their maximum and minimum values, with $e_k = \alpha_E^{(k)}$.
\end{Theorem}
\begin{proof}
  To show that $\gamma_E^{(j)}$ take alternating values, it is sufficient to show that \[
    \text{sign}((e_j - e_j^{*})(e_{j-1} - 
    e_{j-1}^{*})) = -1
  \] for all $1 \le j \le k$. 
  We know from Lemma \ref{Wn0}  that $W_j(\xi) 
  \neq 0$ for all $j \le k$.
  We start again from 
         \[
    \vphantom{V}^{E}V_{k+1}^{j}(\xi) = 
    (e_j - e_j^{*}) V_j(\xi) W_{k-j}(\xi).
  \] Under our assumptions all the terms on the right-hand side are 
  non-zero. Furthermore we know that $\xi$ 
is smaller than any root of 
$\vphantom{V}^{E}V_{k+1}^{j}$ and so $
\text{sign} (\vphantom{V}^{E}V_{k+1}^{j}
(\xi)) = (-1)^{k+1}$. Thus,
  \[
    \text{sign}(e_j - e_j^{*}) = (-1)^{k+1}
    \text{sign}(
    V_j(\xi) W_{k-j}(\xi)
    ).
  \] 
  Replacing $j$ with  $j-1$ we also obtain
  \[
    \text{sign}(e_{j-1} - e_{j-1}^{*}) = 
    (-1)^{k+1}
    \text{sign}(
    V_{j-1}(\xi) W_{k+1-j}(\xi)
    ).
  \] Putting these two results together, we obtain 
  \begin{equation}
    \label{abc_1}
    \text{sign}((e_j- e_j^{*})(e_{j-1} - 
    e_{j-1}^{*})) = \text{sign}(
    V_{j}(\xi) V_{j-1}(\xi) W_{k-j}(\xi) W_{k+1-j}(\xi)
    ).
  \end{equation}
  Now we evaluate \eqref{dec} at $\xi$ and 
  obtain 
  \[
    V_j(\xi)W_{k+1-j}(\xi) =
    (\xi - 1) \gamma_R^{(j)}
    V_{j-1}(\xi)W_{k-j}(\xi).
  \]
 Since under our assumptions all terms 
are non-zero, we obtain 
\begin{equation}
  \label{abc_2}  
  \text{sign}( V_j(\xi) V_{j-1}(\xi) W_{k-j}(\xi) W_{k+1-j}(\xi) ) = 
  \text{sign}(\xi - 1) = -1.
\end{equation}
 Combining \eqref{abc_1} and \eqref{abc_2}, we obtain the result. 

To show that $e_k = \alpha_E^{(k)}$, we evaluate \eqref{EV} at $\xi$ and obtain  \[
\vphantom{V}^{E}V^{k}_{k+1}(\xi) =
(e_k - e_k^{*})
V_k(\xi). \]
Since $\xi$ lies outside of $\mathcal{I}_k$ it follows that $V_k$ has no root less than or equal than $\xi$. Therefore, $\text{sign}(V_k(\xi)) = (-1)^{k}$. 
Similarly, we see that $\text{sign}(\vphantom{V}^{E}V^{k}_{k+1}(\xi) ) = (-1)^{k+1}$. Thus \[
\text{sign}(e_k - e_k^{*}) = (-1)^{k+1}(-1)^{k} = -1,
\] and so $e_k = \alpha_E^{(k)}$.
\end{proof}

\subsubsection{Choice of $\gamma_R^{(j)}$}

\begin{Theorem}
\label{choicer}
Let $e_j$, $r_j$ be parameters that give the smallest possible real root $\xi$ of $V_{k+1}$. Suppose that $\xi$ lies outside of $\mathcal{I}_{\leq k} $ and $\xi < 1$. 
Then this is achieved by letting the parameters $\gamma_R^{(j)}$ alternate between their maximum and minimum values, with $r_k = \alpha_R^{(k)}$.
\end{Theorem}
\begin{proof}
It is sufficient to prove that \[
\text{sign}(r_j - r_j^{*}) = \text{sign}(e_j - e_j^{*})
\] for all $ 1 \leq j \leq k$. 
In the proof of Theorem \ref{choice} it has been shown that \[
\text{sign}(e_j - e_j^{*}) = (-1)^{k+1}
\text{sign}(
V_j(\xi) W_{k-j}(\xi)
).
\]
The same argument yields that \[
\text{sign}(r_j^{*} - r_j) = (-1)^{k+1}
\text{sign}(
V_{j-1}(\xi) W_{k-j}(\xi)
).
\] From these observations, we obtain \[
\text{sign}(r_j^{*} - r_j) = \text{sign}(e_j - e_j^{*})\text{sign}(V_j(\xi)V_{j-1}(\xi)).
\] Since $\xi$ lies outside of $\mathcal{I}_{\leq k} $ it follows that  $\text{sign}(V_l(\xi)) =
(-1)^{l}$ for all $l \le k$. Thus \[
\text{sign}(r_j^{*} - r_j) = \text{sign}(e_j - e_j^{*})(-1)^{j}(-1)^{j-1} = \text{sign}
(e_j^{*} - e_j)
.\]
\end{proof}

\begin{Remark}
Based on the previous results we may compute the left endpoint of
$\mathcal{I}_{k+1} = [\min\{ a_{k+1}, 1 \}, \max\{ b_{k+1}, 1 \}]$, recursively as  $a_{k+1} = \min\{a_k, \xi\}$, where
$\xi$ is the smallest real root of
$V_{k+1}$ for the choice of the parameters 
described in Theorems \ref{choice} and \ref{choicer}, i.e., alternating between their minimum and maximum values, ending in all cases with  $\gamma_E^{(k)} = \alpha_E^{(k)}$ and
$\gamma_R^{(k)} = \alpha_R^{(k)}$. The right endpoint $b_{k+1}$ is simply the largest real root of $V_{k+1}$ where all parameters take their maximum value.
\end{Remark}


\section{Numerical experiments}
\label{sec:NumExpts}

\subsection{Randomly-generated matrices}
We now undertake numerical tests to validate the theoretical bounds developed in the previous sections, using {\scshape Matlab} R2022a. We first determine the extremal eigenvalues of $\mathcal{P}^{-1}\mathcal{A}$ on randomly-generated linear systems. Specifically, we consider a simplified case with $A_i = E_i \equiv 0$, $i \ge 1$, and run 
a number of different test cases, involving a number of combinations of values for the extremal eigenvalues of the symmetric positive definite matrices involved, which are stated in Table \ref{tab:1}. 
In particular,
\begin{itemize}
	\item For $N = 2$, we have 3 parameters with 9 combinations for the different endpoints of the corresponding intervals. Overall, we run $9^3=729$ test cases.
	\item For $N = 3$, we have 4 parameters with 4 combinations for the different endpoints of the corresponding intervals. Overall, we run $4^4=256$ test cases.
	\item For $N = 4$, we have 5 parameters with 4 combinations for the different endpoints of the corresponding intervals. Overall, we run $4^5=1024$ test cases. 
\end{itemize}

\begin{table}[h!]
    \caption{Extremal eigenvalues of the relevant symmetric positive definite matrices used in the verification of the bounds.}
    \label{tab:1}
    \centering
    \begin{tabular}{l|ccc||c|cc}
    $\alpha_E^{(0)}$, $\alpha_R^{(1)}$, $\alpha_R^{(2)}$ & 0.1 & 0.3 & 0.9 &
    $\alpha_E^{(0)}$, $\alpha_R^{(1)}$, \ldots, $\alpha_R^{(N)}$ & 0.1 & 0.9\\
    $\beta_E^{(0)}$, $\beta_R^{(1)}$, $\beta_R^{(2)}$ & 1.2 & 1.8 & 5 &
    $\beta_E^{(0)}$, $\beta_R^{(1)}$, \ldots, $\beta_R^{(N)}$ & 1.2 & 5\\
    \hline
    \multicolumn{4}{c||}{Case $N = 2$} & \multicolumn{3}{c}{Cases $N > 2$} \\
    \end{tabular}
\end{table}
We run each test case $10$ times, generating random matrices that satisfy the relevant spectral properties, and we record the most extreme eigenvalues for each test case. 
In more detail, the dimensions $n_0$, $n_1, \ldots, n_N$ are computed using the closest integer to $50+10\cdot\texttt{rand}$, using {\scshape Matlab}'s \texttt{rand} function, recomputing as necessary to ensure that $n_k \leq n_{k-1}$, $k = 1,\ldots, N$. 
The matrices $A_k$ and $B_k$ are computed using {\scshape Matlab}'s \texttt{randn} function, whereupon we take the symmetric part of $A_k$ and then add an identity matrix multiplied by the largest absolute value of the negative eigenvalues (multiplied by $1.01$ if $k=0$), to ensure symmetric positive semi-definiteness (definiteness if $k=0$). We then choose $\widehat{S}_0$ as a linear combination of $A$ and the identity matrix, such that the eigenvalues of $E_0$ are contained in $[\alpha_E^{(0)}, \beta_E^{(0)}]$, following the procedure in \cite[Sec. 3]{pearson2023symmetric}. We construct $\widehat{S}_k$ for $k \ge 1$ similarly, obtaining that the eigenvalues of $R_k R_k ^\top
= \widehat S_k^{-1/2} B_k \widehat S_{k-1}^{-1} B_k^\top S_k^{-1/2}$ are in 
the interval $[\alpha_R^{(k)}, \beta_R^{(k)}]$. 

%

\subsection*{Multiple saddle-point systems with $N=2$.}
\begin{figure}[h!]
    \centering
    \includegraphics[width=.48\linewidth]{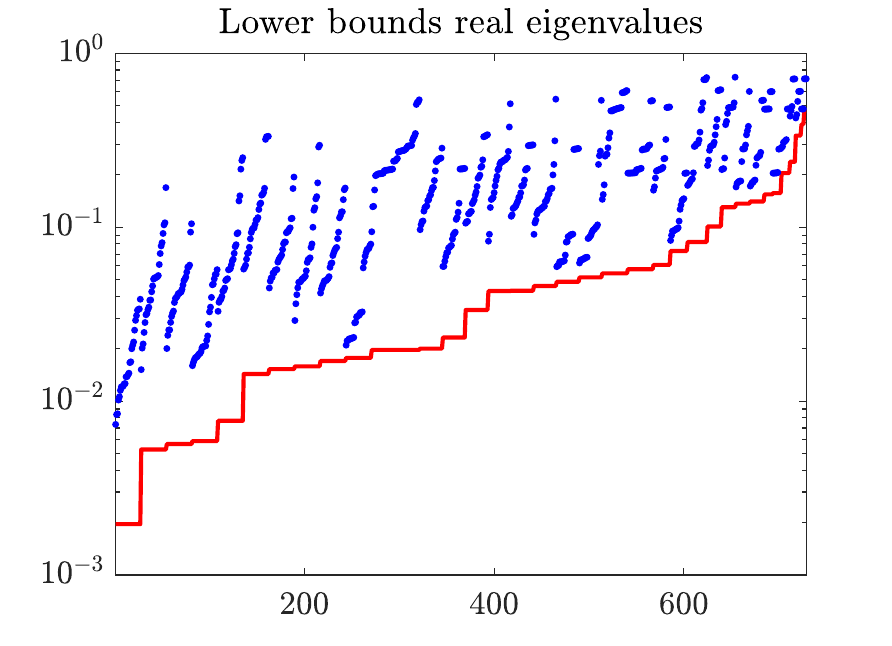}
    \includegraphics[width=.48\linewidth]{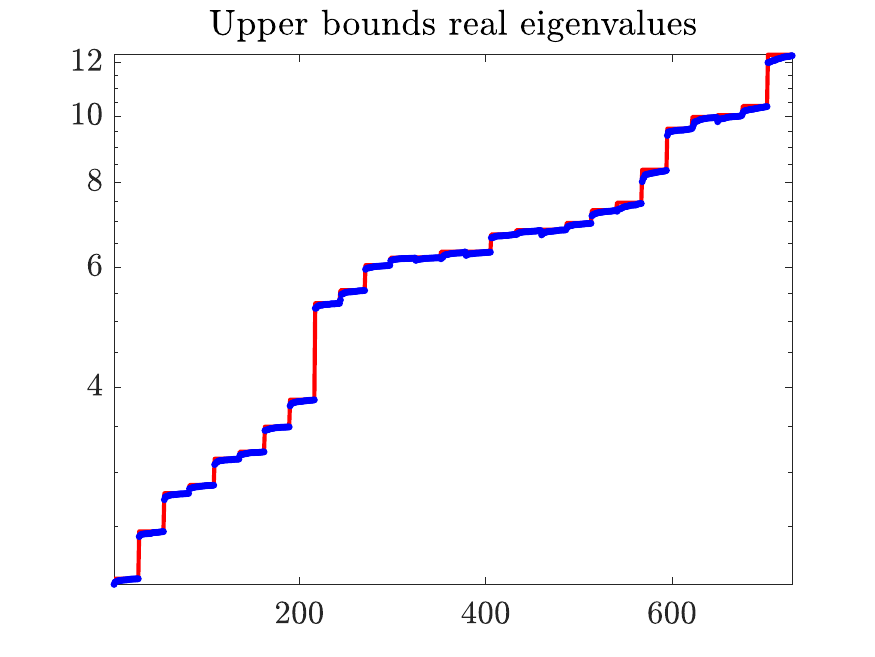}
        \caption{Multiple saddle-point linear system with $N = 2$ (three blocks). Extremal eigenvalues of the preconditioned matrix (blue dots) and bounds (red line) after $10$ runs with each combination of the parameters from Table \ref{tab:1}.} 
    \label{fig:test1}
\end{figure}

In Figure \ref{fig:test1} we report the extremal real eigenvalues of $\mathcal{P}^{-1} \mathcal{A}$, compared to the theoretical bounds of Section \ref{sec:roots}, for the double saddle-point setting $N=2$.
The extremal eigenvalues are sorted (and the computed bounds accordingly) for improved readability. 
We notice
that the plots indicate (for these problems) that the bounds capture
the behaviour of the eigenvalues very well, particularly so for the largest eigenvalues, which are almost perfectly described
by our estimates.
Regarding the complex eigenvalues, we computed for each of the 729 test cases the largest value of $|\lambda-1|$, i.e., in the complex plane the radius of the circle of center
$(1,0)$ containing all complex eigenvalues. We obtained values ranging from $0.197$ to $0.988$.

\subsection*{Multiple saddle-point systems with $N=3$.}
\begin{figure}[h!]
    \centering
    \includegraphics[width=.48\linewidth]{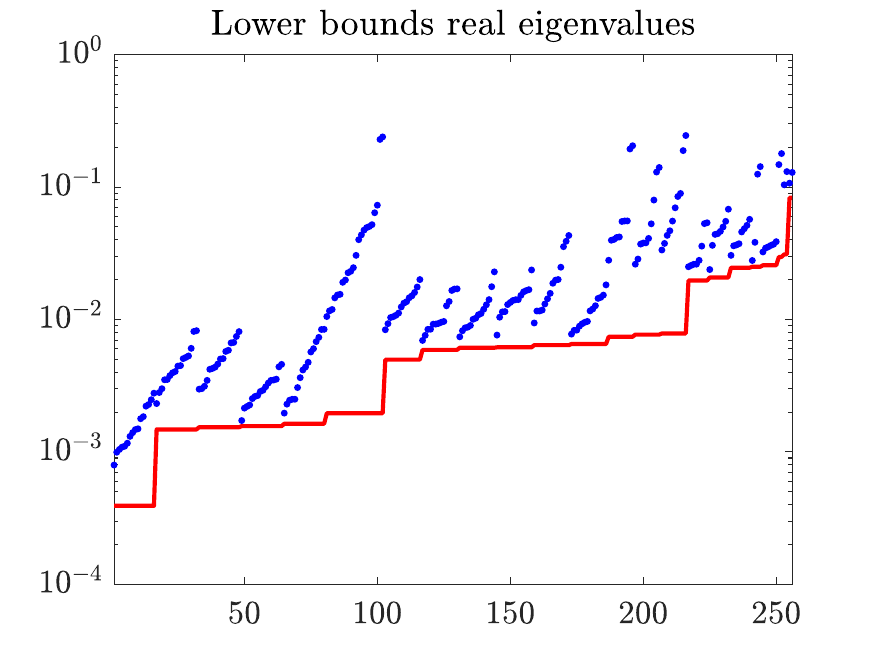}
    \includegraphics[width=.48\linewidth]{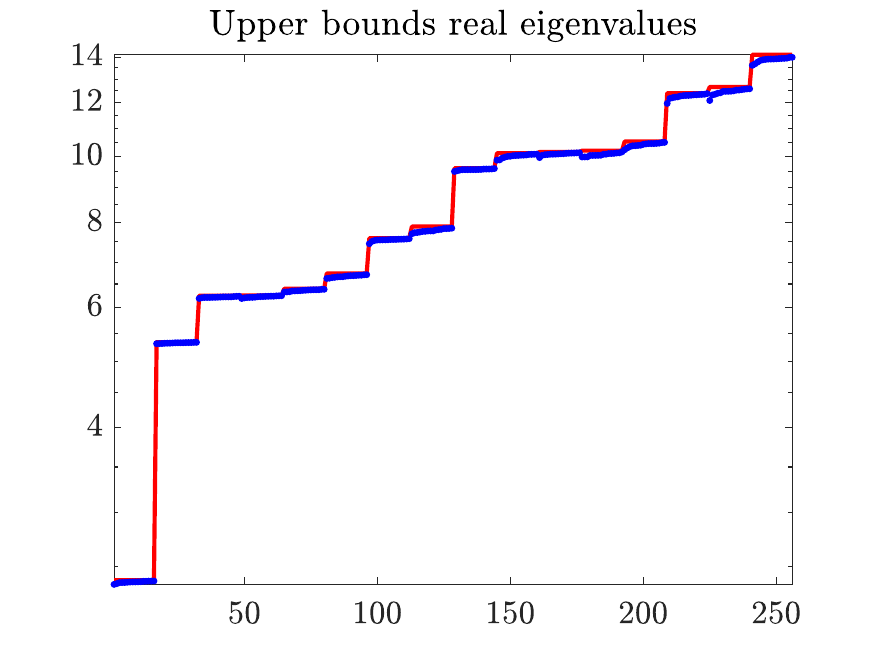}
        \caption{Multiple saddle-point linear system with $N = 3$ (four blocks). Extremal eigenvalues of the preconditioned matrix (blue dots) and bounds (red line) after $10$ runs with each combination of the parameters from Table \ref{tab:1}.}
    \label{fig:test2}
\end{figure}

In Figure \ref{fig:test2} we report the extremal real eigenvalues of $\mathcal{P}^{-1} \mathcal{A}$, compared to the theoretical bounds of Section \ref{sec:roots}, for the triple saddle-point setting $N=3$. 
The theory again captures the practical behaviour very well, in particular for the largest eigenvalues.
For the complex eigenvalues, the largest value of $|\lambda-1|$ for each of the 256 test cases ranged from $0.675$ to $0.99994$.
\

\newpage
\subsection*{Multiple saddle-point systems with $N=4$.}
\begin{figure}[h!]
    \centering
    \includegraphics[width=.48\linewidth]{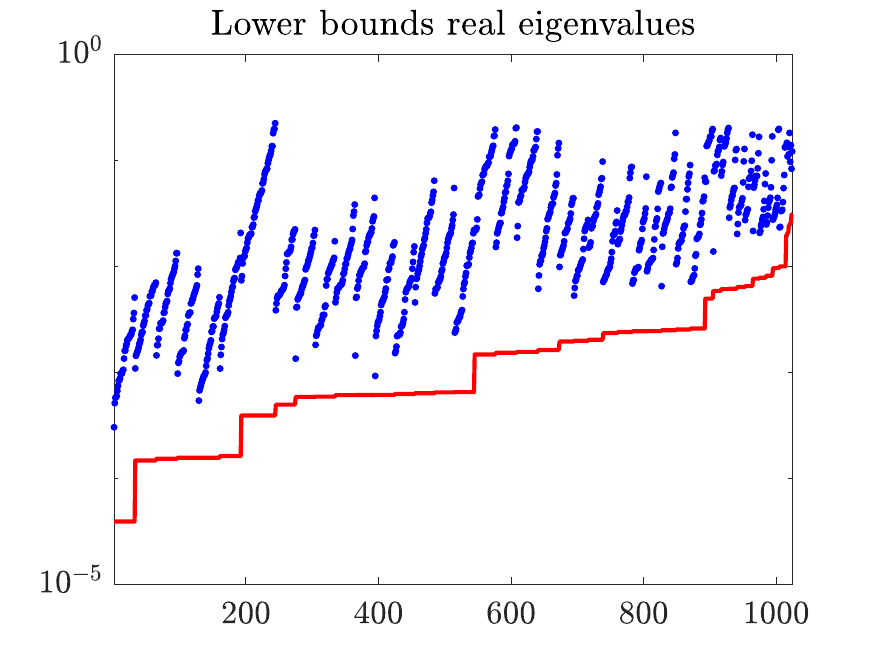}
    \includegraphics[width=.48\linewidth]{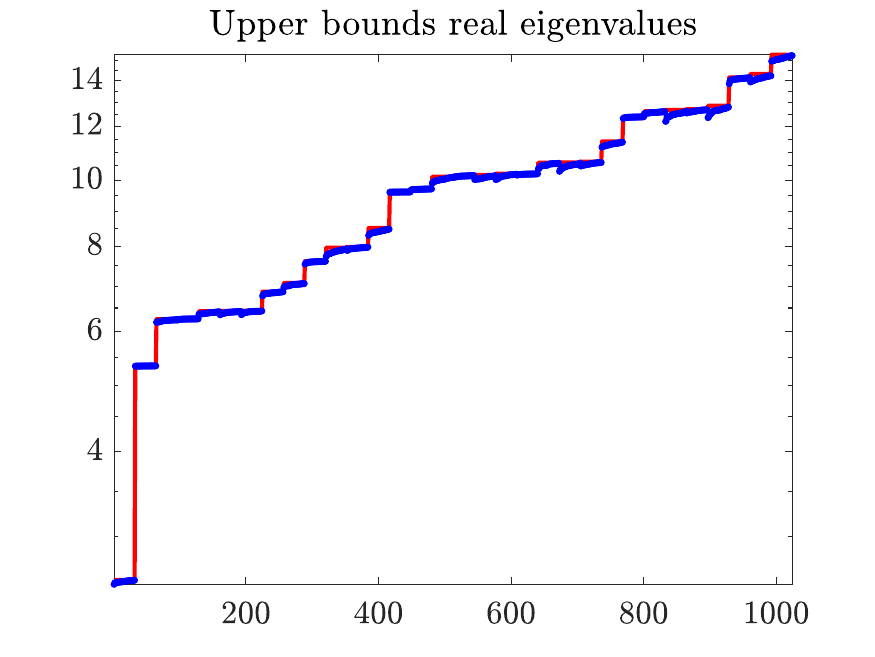}
	\caption{Multiple saddle-point linear system with $N = 4$ (five blocks). Extremal eigenvalues of the preconditioned matrix (blue dots) and bounds (red line) after $10$ runs with each combination of the parameters from Table \ref{tab:1}.}
    \label{fig:test3}
\end{figure}

Figure \ref{fig:test3} shows the extremal eigenvalues of $\mathcal{P}^{-1} \mathcal{A}$, compared to the theoretical bounds of Section \ref{sec:roots}, for the case $N=4$. 
We see that our theoretical bounds describe the behaviour of the real eigenvalues well, particularly for the largest eigenvalues, and for the 1024 test cases the largest value of $|\lambda-1|$ for the complex eigenvalues ranged from {$0.750$} to $0.9999985$. 

\subsection*{Role of the real eigenvalues on GMRES convergence.} 
\begin{figure}[h!]
    \centering
    \includegraphics[width=.335\linewidth]{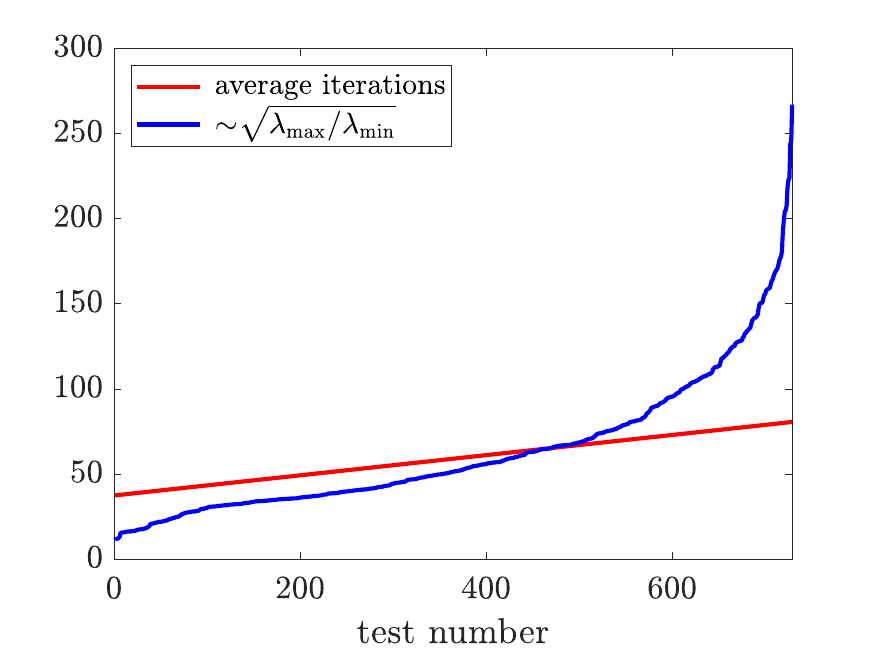}
	\hspace{-3mm}
    \includegraphics[width=.335\linewidth]{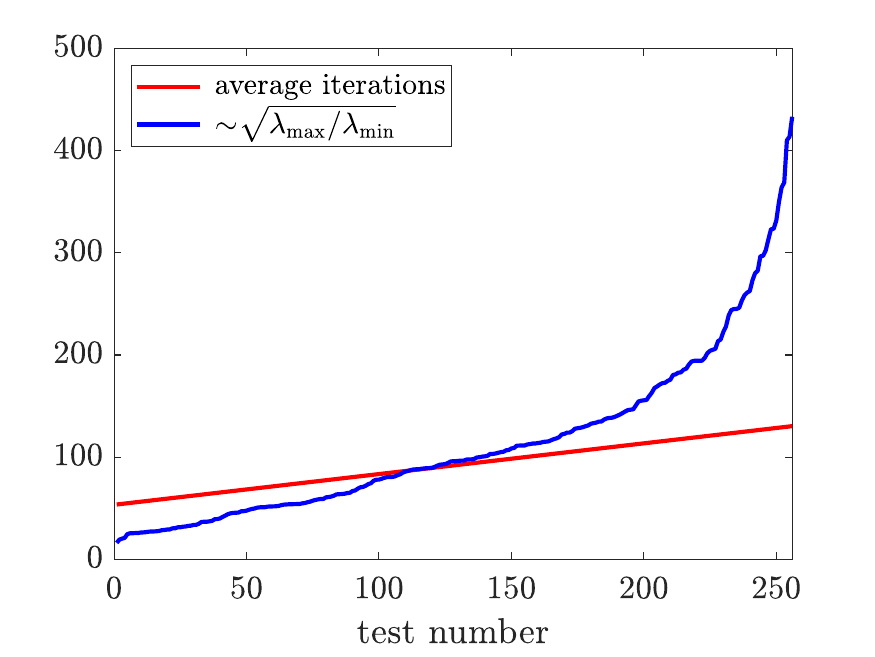}
	\hspace{-3mm}
    \includegraphics[width=.335\linewidth]{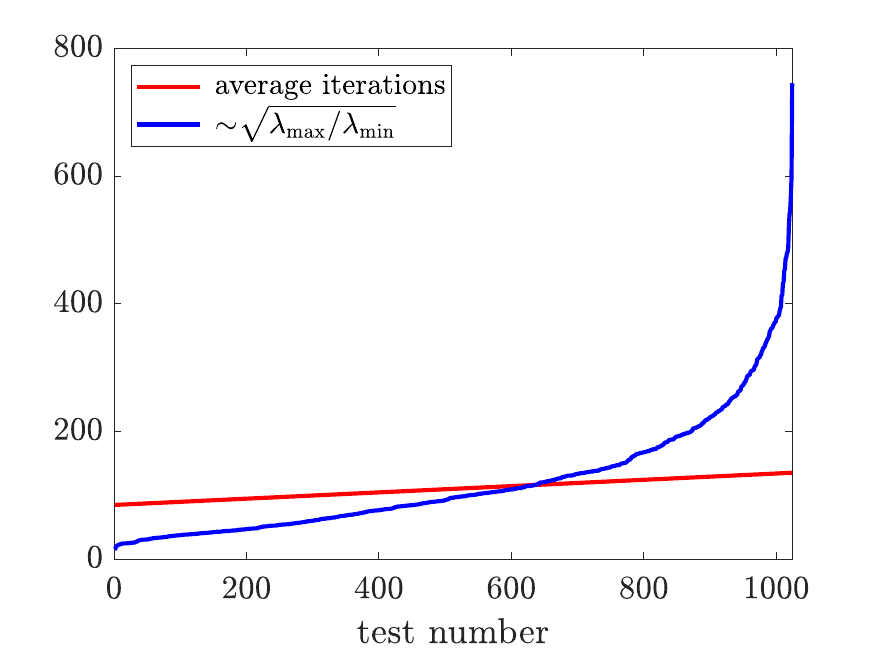}
	\caption{Plot of average number of GMRES iterations and (scaled) square root of the ratio between the largest and the smallest real eigenvalues, sorted across all	test case numbers. From left to right: $N = 2$, $N = 3$, and $N = 4$.}
	\label{its_vs_cond}
\end{figure}

To further analyze the role of real and complex eigenvalues in the efficacy of the block-triangular preconditioner we also run the GMRES method (without restarting) for all 
test cases. 
For the purposes of this exercise, the right-hand side is pre-computed such that the exact solution is the vector of ones. For each test, we record the average number of iterations (to obtain the relative residual norm less than $10^{-10}$) over
the 10 runs, and the following indicator, based on real eigenvalues only:
\[ 
\frac{\max\{\lambda: \lambda \in \R\}}{\min\{ \lambda: \lambda \in \R\}} =: \frac{\lambda_{\max}}{\lambda_{\min}}. \]
In Figure \ref{its_vs_cond}, we plot both the (uniformly-scaled) values of $\sqrt \frac{\lambda_{\max}}{\lambda_{\min}}$ for each problem, and a linear fit for the numbers of GMRES iterations required for these problems, for multiple
saddle-point linear systems with $N = 2, 3, 4$.
From these plots we deduce that, as may be expected, the ratio between the largest and the smallest real eigenvalue of the preconditioned matrix appears to influence
the numbers of GMRES iterations.

\subsection{PDE-constrained optimization example}

We now test the results of Section \ref{sec:roots} on the following PDE-constrained optimization example (see work such as \cite[Sec. 5]{BMPP_COAP25}, \cite{MNN}, and \cite[Sec. 4.1]{pearson2023symmetric}):
\begin{equation*}
\min_{y,u} ~~ \frac{1}{2}\left\|y-\widehat{y}\right\|_{L^2(\partial\Omega)}^2+\frac{\beta}{2}\left\|u\right\|_{L^2(\Omega)}^2 \qquad \text{s.t.} \quad \left\{\begin{array}{rl}
-\Delta y+y+u=0 & \text{in }\Omega, \\
\frac{\partial y}{\partial n}=0 & \text{on }\partial\Omega. \\
\end{array}\right.
\end{equation*}
Here, the domain $\Omega := (0,1)^2$ with boundary $\partial\Omega$, $y$ and $u$ denote \emph{state} and \emph{control variables}, and $\beta$ is a (positive) regularization parameter. The function $\widehat{y}$ denotes a (given) \emph{desired state}, which for this example is obtained by solving the forward problem with control $4x_1(1-x_1)+x_2$, with $x_1$ and $x_2$ denoting the spatial coordinates.

Discretizing this problem using $P_1$ finite element basis functions yields the double saddle-point linear system
\begin{equation}\label{PDECO_system}
\left[\begin{array}{ccc}
\beta M & M & 0 \\ M & 0 & L \\ 0 & L & M_{\partial\Omega} \\
\end{array}\right]\left[\begin{array}{c}
u_h \\ p_h \\ y_h \\
\end{array}\right]=\left[\begin{array}{c}
0 \\ 0 \\ \widehat{y}_h \\
\end{array}\right],
\end{equation}
where $y_h$, $u_h$, $p_h$, and $\widehat{y}_h$ denote the discretized state, control, \emph{adjoint}, and desired state, $M$ and $M_{\partial\Omega}$ denote a finite element mass matrix and boundary mass matrix, and $L$ is the sum of finite element stiffness and mass matrices. In this experiment, we precondition \eqref{PDECO_system} with
\begin{equation*}
\mathcal{P} = \left[\begin{array}{ccc}
\beta \widehat{M} & M & 0 \\ 0 & -\frac{1}{\beta} \widehat{M} & L \\ 0 & 0 & \beta \widehat{L} M^{-1} \widehat{L} \\
\end{array}\right],
\end{equation*}
where $\widehat{M}$ denotes a number (\texttt{Cheb}) of Chebyshev semi-iterations applied to $M$ (see \cite{GVI,GVII,WathenRees}) and $\widehat{L}$ denotes two V-cycles of the \texttt{HSL\_MI20} algebraic multigrid routine applied to $L$ (see \cite{HSL_MI20,HSL_MI20_code}). Our theoretical bounds may be obtained using known results \cite[Sec. 5]{BMPP_COAP25}. Specifically,
\begin{equation*}
[\alpha_E^{(0)},\beta_E^{(0)}] \in [1-\omega, 1+\omega], \qquad \text{where} \quad \omega = 1\left/T_{\texttt{Cheb}} \left( \frac{5}{3} \right),\right.
\end{equation*}
with $T_{\texttt{Cheb}}$ the $\texttt{Cheb}$-th Chebyshev polynomial, $\alpha_R^{(1)} = (\alpha_E^{(0)})^2$, $\beta_R^{(1)} = (\beta_E^{(0)})^2$, $\alpha_E^{(1)} = \beta_E^{(1)} = \alpha_E^{(2)} \equiv 0$. We compute the values of $\beta_E^{(2)}$, $\alpha_R^{(2)}$, and $\beta_R^{(2)}$ explicitly.

\begin{table}[h!]
\centering
\caption{Computed extremal real eigenvalues of preconditioned system (denoted `\text{Comp}') and theoretical bounds (denoted `\text{Bound}'), for PDE-constrained optimization problem with $h = 2^{-5}$, $\beta \in \{ 1, 10^{-3} \}$, and a range of Chebyshev semi-iterations \texttt{Cheb}. Here, $\lambda_c$ denotes the computed value of $\max_{\lambda,\lambda\notin\mathbb{R}}\,|\lambda-1|$.}\label{Results_Eigenvalues_h2-5}
\begin{tabular}{|c||c|c||c|c||c|}
\hline
 & \multicolumn{5}{c|}{$\beta = 1$} \\ \cline{2-6}
\texttt{Cheb} & $\text{Bound}_l$ & $\text{Comp}_l$ & $\text{Comp}_u$ & $\text{Bound}_u$ & $\lambda_c$ \\ \hline \hline
1 & 0.0544 & 0.3082 & 6.4197 & 6.4218 & 0.9313 \\ \hline
3 & 0.5216 & 0.6793 & 5.0992 & 5.1057 & 0.4809 \\ \hline
5 & 0.7769 & 0.8263 & 4.9910 & 5.0073 & 0.2196 \\ \hline
10 & 0.9502 & 0.9627 & 4.9780 & 4.9963 & 0.0347 \\ \hline
\end{tabular}

\begin{tabular}{|c||c|c||c|c||c|}
\hline
 & \multicolumn{5}{c|}{$\beta = 10^{-3}$} \\ \cline{2-6}
\texttt{Cheb} & $\text{Bound}_l$ & $\text{Comp}_l$ & $\text{Comp}_u$ & $\text{Bound}_u$ & $\lambda_c$ \\ \hline \hline
1 & 0.0544 & 0.3082 & 3984.4 & 3984.4 & 0.9313 \\ \hline
3 & 0.5216 & 0.6797 & 3983.9 & 3983.9 & 0.4809 \\ \hline
5 & 0.7769 & 0.8272 & 3983.8 & 3983.8 & 0.2196 \\ \hline
10 & 0.9502 & 0.9636 & 3983.8 & 3983.8 & 0.0329 \\ \hline
\end{tabular}
\end{table}

\begin{table}
\centering
\caption{Computed extremal real eigenvalues of preconditioned system (denoted `\text{Comp}') and theoretical bounds (denoted `\text{Bound}'), for PDE-constrained optimization problem with $h = 2^{-6}$, $\beta \in \{ 1, 10^{-3} \}$, and a range of Chebyshev semi-iterations \texttt{Cheb}. Here, $\lambda_c$ denotes the computed value of $\max_{\lambda,\lambda\notin\mathbb{R}}\,|\lambda-1|$.}\label{Results_Eigenvalues_h2-6}
\begin{tabular}{|c||c|c||c|c||c|}
\hline
 & \multicolumn{5}{c|}{$\beta = 1$} \\ \cline{2-6}
\texttt{Cheb} & $\text{Bound}_l$ & $\text{Comp}_l$ & $\text{Comp}_u$ & $\text{Bound}_u$ & $\lambda_c$ \\ \hline \hline
1 & 0.0544 & 0.3082 & 6.4009 & 6.4082 & 0.9313 \\ \hline
3 & 0.5214 & 0.6791 & 5.0803 & 5.0907 & 0.4809 \\ \hline
5 & 0.7758 & 0.8260 & 4.9722 & 4.9935 & 0.2196 \\ \hline
10 & 0.9480 & 0.9577 & 4.9593 & 4.9820 & 0.0373 \\ \hline
\end{tabular}

\begin{tabular}{|c||c|c||c|c||c|}
\hline
 & \multicolumn{5}{c|}{$\beta = 10^{-3}$} \\ \cline{2-6}
\texttt{Cheb} & $\text{Bound}_l$ & $\text{Comp}_l$ & $\text{Comp}_u$ & $\text{Bound}_u$ & $\lambda_c$ \\ \hline \hline
1 & 0.0544 & 0.3082 & 3969.3 & 3969.3 & 0.9313 \\ \hline
3 & 0.5124 & 0.6791 & 3968.8 & 3968.8 & 0.4809 \\ \hline
5 & 0.7758 & 0.8263 & 3968.7 & 3968.7 & 0.2196 \\ \hline
10 & 0.9480 & 0.9619 & 3968.7 & 3968.7 & 0.0333 \\ \hline
\end{tabular}
\end{table}

\begin{table}
\centering
\caption{Number of GMRES iterations required for solution of PDE-constrained optimization problem, for a range of Chebyshev semi-iterations \texttt{Cheb}.}\label{ResultsTable_Iterations}
\begin{tabular}{|c||c|c|c|c|}
\hline
$h$ & $2^{-5}$ & $2^{-5}$ & $2^{-6}$ & $2^{-6}$ \\ \hline
$\texttt{Cheb} \backslash \beta$ & $~1~$ & $~10^{-3}~$ & $~1~$ & $~10^{-3}~$ \\ \hline \hline
1 & 50 & 33 & 55 & 37 \\ \hline
3 & 17 & 17 & 16 & 17 \\ \hline
5 & 11 & 14 & 11 & 14 \\ \hline
10 & 6 & 13 & 6 & 13 \\ \hline
\end{tabular}
\end{table}

In Table \ref{Results_Eigenvalues_h2-5} we show the smallest and largest computed real eigenvalues ($\text{Comp}_l$ and $\text{Comp}_u$) of $\mathcal{P}^{-1} \mathcal{A}$ with mesh parameter $h = 2^{-5}$ and $\beta \in \{ 1, 10^{-3} \}$, with different numbers of Chebyshev semi-iterations applied to approximate $M^{-1}$, as well as the analytic bounds ($\text{Bound}_l$ and $\text{Bound}_u$) obtained using the above theory. Table \ref{Results_Eigenvalues_h2-6} displays the analogous results for $h = 2^{-6}$. Table \ref{ResultsTable_Iterations} shows the GMRES iteration numbers required to solve all systems tested, to a relative tolerance of $10^{-10}$. The linear systems with $h = 2^{-5}$ and $h = 2^{-6}$ have dimensions of 3267 and 12675, respectively. Tables \ref{Results_Eigenvalues_h2-5} and \ref{Results_Eigenvalues_h2-6} demonstrate the descriptiveness of the analytic bounds derived in this work on the real eigenvalues of the preconditioned double saddle-point system, for a practical problem. The bounds on the real eigenvalues are tighter for larger values of \texttt{Cheb}, as is the closeness of the complex eigenvalues to $1$. As expected, the more accurate approximation of the mass matrices arising from increasing \texttt{Cheb} also leads to the reduction of GMRES iteration numbers seen in Table \ref{ResultsTable_Iterations}. For larger $\beta$, this reduction is particularly drastic, while for smaller $\beta$ the second Schur complement $S_2$ is less well approximated, leading to larger values of the maximum real eigenvalues and iterations required.

\section{Conclusions}
\label{sec:conc}

In this work, we have categorized the behaviour of complex and real eigenvalues of multiple saddle-point systems following application of a classical block-triangular preconditioner, involving approximations of the top-left block and all successive Schur complements. The bounds on the complex eigenvalues are expressed in terms of the distance from the value $1$ in the complex plane; these involve properties of the leading blocks of the corresponding eigenvectors, but in certain settings may be used to provide practical bounds regardless of the eigenvector structure. The real eigenvalues are precisely characterized in terms of roots of a defined sequence of polynomials, and we have justified how the individual coefficients may be selected in order to achieve the extremal values of the roots, which may be computed very cheaply.

Numerical experiments demonstrate that the theoretical bounds effectively characterize the eigenvalues of the preconditioned system for various numbers of blocks; the bounds are particularly effective when the individual Schur complements are reasonably well approximated, and illustrate that good approximations of each block lead to clustered eigenvalues of the preconditioned system. Tests on both randomly-generated problems and PDE-constrained optimization examples show that in practice this leads to reliable and fast convergence of GMRES, and that the bound on the maximum real eigenvalue, in particular, captures the practical behaviour very closely.

\section*{Acknowledgements}

JWP acknowledges funding from the UK’s Engineering and Physical
Sciences Research Council (EPSRC grant EP/Z533786/1). The authors acknowledge funding from the University of Padua supporting a visit of JWP to the University, during which some of this work was undertaken. 
LB is member of the \textit{Gruppo Nazionale per il Calcolo Scientifico-Istituto Nazionale di Alta Matematica} (GNCS-INdAM).


\begin{thebibliography}{10}

\bibitem{Balani-et-al-2023b}
{\sc F.~Bakrani~Balani, L.~Bergamaschi, {\'A}.~Mart\'{i}nez, and M.~Hajarian},
  {\em Some preconditioning techniques for a class of double saddle point
  problems}, Numerical Linear Algebra with Applications, 31 (2024), p.~e2551.

\bibitem{Balani-et-al-2023a}
{\sc F.~Bakrani~Balani, M.~Hajarian, and L.~Bergamaschi}, {\em Two block
  preconditioners for a class of double saddle point linear systems}, Applied
  Numerical Mathematics, 190 (2023), pp.~155--167.

\bibitem{BSZ2020}
{\sc A.~Beigl, J.~Sogn, and W.~Zulehner}, {\em Robust preconditioners for
  multiple saddle point problems and applications to optimal control problems},
  SIAM Journal on Matrix Analysis and Applications, 41 (2020), pp.~1590--1615.

\bibitem{Benzi2018}
{\sc F.~P.~A. Beik and M.~Benzi}, {\em Iterative methods for double saddle
  point systems}, SIAM Journal on Matrix Analysis and Applications, 39 (2018),
  pp.~902--921.

\bibitem{BeikBenzi2022}
\leavevmode\vrule height 2pt depth -1.6pt width 23pt, {\em Preconditioning
  techniques for the coupled {S}tokes--{D}arcy problem: spectral and
  field-of-values analysis}, Numerische Mathematik, 150 (2022), pp.~257--298.

\bibitem{BEIK2024403}
{\sc F.~P.~A. Beik, C.~Greif, and M.~Trummer}, {\em On the invertibility of
  matrices with a double saddle-point structure}, Linear Algebra and its
  Applications, 699 (2024), pp.~403--420.

\bibitem{Benzi2026}
{\sc M.~Benzi, M.~Feder, L.~Heltai, and F.~Mugnaioni}, {\em Scalable augmented
  {L}agrangian preconditioners for fictitious domain problems}, Computer
  Methods in Applied Mechanics and Engineering, 450 (2026), p.~Art. 118522.

\bibitem{BenziGolubLiesen2005}
{\sc M.~Benzi, G.~H. Golub, and J.~Liesen}, {\em Numerical solution of saddle
  point problems}, Acta Numerica, 14 (2005), pp.~1--137.

\bibitem{BergaNLAA10}
{\sc L.~Bergamaschi}, {\em On eigenvalue distribution of
  constraint-preconditioned symmetric saddle point matrices}, Numerical Linear
  Algebra with Applications, {19} (2012), pp.~{754--772}.

\bibitem{BB2026}
{\sc L.~Bergamaschi and M.~Bergamaschi}, {\em Eigenvalue bounds for symmetric,
  multiple saddle-point matrices with {SPD} preconditioners}, Algorithms, {19}
  (2026), p.~{Art. 359}.

\bibitem{BerFerMar26}
{\sc L.~Bergamaschi, M.~Ferronato, and {\'A}.~Mart\'{\i}nez}, {\em {Block
  triangular preconditioners for double saddle-point linear systems arising in
  the mixed form of poroelasticity equations}}, SIAM Journal on Matrix Analysis
  and Applications, 47 (2026), pp.~132--157.

\bibitem{BMPP_COAP25}
{\sc L.~Bergamaschi, {\'A}.~Mart\'{i}nez, J.~W. Pearson, and A.~Potschka}, {\em
  Spectral analysis of block preconditioners for double saddle-point linear
  systems with application to {PDE}-constrained optimization}, Computational
  Optimization with Applications, 91 (2025), pp.~423--455.

\bibitem{bergamaschi2026eigenvalue}
\leavevmode\vrule height 2pt depth -1.6pt width 23pt, {\em Eigenvalue bounds
  for preconditioned symmetric multiple saddle-point matrices}, Linear Algebra
  and its Applications,  (2026).

\bibitem{HSL_MI20}
{\sc J.~Boyle, M.~Mihajlovi\'{c}, and J.~Scott}, {\em {HSL\_MI20}: {An}
  efficient {AMG} preconditioner for finite element problems in {3D}},
  International Journal for Numerical Methods in Engineering, 82 (2010),
  pp.~64--98.

\bibitem{Bradley}
{\sc S.~Bradley and C.~Greif}, {\em Eigenvalue bounds for double saddle-point
  systems}, IMA Journal of Numerical Analysis, 43 (2023), pp.~3564--3592.

\bibitem{Mardal2026}
{\sc M.~Cai, M.~Kuchta, J.~Li, Z.~Li, and K.-A. Mardal}, {\em Parameter-robust
  preconditioners for a four-field thermo-poroelasticity model}, SIAM Journal
  on Scientific Computing, 48 (2026), pp.~A49--A73.

\bibitem{Szyld}
{\sc P.~Chidyagwai, S.~Ladenheim, and D.~B. Szyld}, {\em Constraint
  preconditioning for the coupled {S}tokes--{D}arcy system}, SIAM Journal on
  Scientific Computing, 38 (2016), pp.~A668--A690.

\bibitem{desturler}
{\sc E.~de~Sturler and J.~Liesen}, {\em Block-diagonal and constraint
  preconditioners for nonsymmetric indefinite linear systems. {Part I:
  Theory}}, SIAM Journal on Scientific Computing, 26 (2005), pp.~1598--1619.

\bibitem{ESW2014}
{\sc H.~C. Elman, D.~J. Silvester, and A.~J. Wathen}, {\em Finite Elements and
  Iterative Solvers With Applications in Incompressible Fluid Dynamics}, Oxford
  University Press, Oxford, UK, 2nd~ed., 2014.

\bibitem{FerFraJanCasTch19}
{\sc M.~Ferronato, A.~Franceschini, C.~Janna, N.~Castelletto, and H.~A.
  Tchelepi}, {\em {A general preconditioning framework for coupled
  multi-physics problems}}, Journal of Computational Physics, 398 (2019),
  p.~Art. 108887.

\bibitem{FRIGO202136}
{\sc M.~Frigo, N.~Castelletto, M.~Ferronato, and J.~A. White}, {\em Efficient
  solvers for hybridized three-field mixed finite element coupled
  poromechanics}, Computers \& Mathematics with Applications, 91 (2021),
  pp.~36--52.

\bibitem{GVII}
{\sc G.~H. Golub and R.~S. Varga}, {\em Chebyshev semi-iterative methods,
  successive over-relaxation iterative methods, and second order {R}ichardson
  iterative methods. {Part I}}, Numerische Mathematik, 3 (1961), pp.~147--156.

\bibitem{GVI}
\leavevmode\vrule height 2pt depth -1.6pt width 23pt, {\em Chebyshev
  semi-iterative methods, successive over-relaxation iterative methods, and
  second order {R}ichardson iterative methods. {Part II}}, Numerische
  Mathematik, 3 (1961), pp.~157--168.

\bibitem{greif2026}
{\sc C.~Greif}, {\em A {BFBt} preconditioner for double saddle-point systems},
  IMA Journal of Numerical Analysis,  (2026).
\newblock draf154.

\bibitem{Greif2017}
{\sc C.~Greif, S.~He, and P.~Liu}, {\em {SYM-ILDL}: {I}ncomplete {LDLT}
  factorization of symmetric indefinite and skew-symmetric matrices}, ACM
  Transactions on Mathematical Software, 44 (2017), p.~Art. 1.

\bibitem{greifhe2023}
{\sc C.~Greif and Y.~He}, {\em Block preconditioners for the marker-and-cell
  discretization of the {S}tokes--{D}arcy equations}, SIAM Journal on Matrix
  Analysis and Applications, 44 (2023), pp.~1540--1565.

\bibitem{HSL_MI20_code}
{\sc {HSL Mathematical Software Library}}, {\em \texttt{HSL\_MI20} unsymmetric
  system: algebraic multigrid preconditioner, code available at
  \url{https://www.hsl.rl.ac.uk/catalogue/hsl_mi20.html}}, 2023.

\bibitem{MNN}
{\sc K.-A. Mardal, B.~F. Nielsen, and M.~Nordaas}, {\em Robust preconditioners
  for {PDE}-constrained optimization with limited observations}, BIT Numerical
  Mathematics, 57 (2017), pp.~405--431.

\bibitem{Mardal2026-MMMAS}
{\sc K.-A. Mardal, J.~Sogn, and S.~Takacs}, {\em A robust and time-parallel
  preconditioner for parabolic reconstruction problems using {I}sogeometric
  analysis}, Mathematical Models and Methods in Applied Sciences, 36 (2026),
  pp.~861--886.

\bibitem{minres}
{\sc C.~C. Paige and M.~A. Saunders}, {\em Solution of sparse indefinite
  systems of linear equations}, SIAM Journal on Numerical Analysis, 12 (1975),
  pp.~617--629.

\bibitem{PPNLAA24}
{\sc J.~W. Pearson and A.~Potschka}, {\em Double saddle-point preconditioning
  for {K}rylov methods in the inexact sequential homotopy method}, Numerical
  Linear Algebra with Applications, 31 (2024), p.~e2553.

\bibitem{pearson2023symmetric}
\leavevmode\vrule height 2pt depth -1.6pt width 23pt, {\em {On symmetric
  positive definite preconditioners for multiple saddle-point systems}}, IMA
  Journal of Numerical Analysis, 44 (2024), pp.~1731--1750.

\bibitem{Pilotto-et-al2026}
{\sc M.~Pilotto, L.~Bergamaschi, and {\'A}.~Mart\'{i}nez}, {\em Spectral
  analysis of block diagonally preconditioned multiple saddle-point matrices
  with inexact {S}chur complements}, Numerical Linear Algebra with
  Applications, 33 (2026), p.~e70110.

\bibitem{RamGar2023}
{\sc A.~Ramage and E.~C. {Gartland Jr.}}, {\em A preconditioned nullspace
  method for liquid crystal director modeling}, SIAM Journal on Scientific
  Computing, 35 (2013), pp.~B226--B247.

\bibitem{Rhebergen}
{\sc S.~Rhebergen, G.~N. Wells, A.~J. Wathen, and R.~F. Katz}, {\em Three-field
  block preconditioners for models of coupled magma/mantle dynamics}, SIAM
  Journal on Scientific Computing, 37 (2015), pp.~A2270--A2294.

\bibitem{saadbook}
{\sc Y.~Saad}, {\em Iterative Methods for Sparse Linear Systems}, SIAM,
  Philadelphia, PA, 2nd~ed., 2003.

\bibitem{saad1986gmres}
{\sc Y.~Saad and M.~H. Schultz}, {\em G{MRES}: a generalized minimal residual
  algorithm for solving nonsymmetric linear systems}, SIAM Journal on
  Scientific and Statistical Computing, 7 (1986), pp.~856--869.

\bibitem{Sim04}
{\sc V.~Simoncini}, {\em Block triangular preconditioners for symmetric
  saddle-point problems}, Applied Numerical Mathematics, 49 (2004), pp.~63--80.

\bibitem{SZ}
{\sc J.~Sogn and W.~Zulehner}, {\em {Schur complement preconditioners for
  multiple saddle point problems of block tridiagonal form with application to
  optimization problems}}, IMA Journal of Numerical Analysis, 39 (2018),
  pp.~1328--1359.

\bibitem{WathenRees}
{\sc A.~J. Wathen and T.~Rees}, {\em Chebyshev semi-iteration in
  preconditioning for problems including the mass matrix}, Electronic
  Transactions on Numerical Analysis, 34 (2009), pp.~125--135.

\end{thebibliography}
\end{document}